\documentclass[ijoc,nonblindrev]{informs3} 

\DoubleSpacedXII

\usepackage{natbib}
 \bibpunct[, ]{(}{)}{,}{a}{}{,}%
 \def\bibfont{\small}%

\usepackage{enumitem}
 \usepackage{multirow}

\TheoremsNumberedThrough     

\EquationsNumberedThrough    

\begin{document}


\RUNAUTHOR{Zhang}

\RUNTITLE{Integer Programming Approaches for Risk-Adjustable DRCCs}

\TITLE{Integer Programming Approaches for
Distributionally Robust Chance Constraints
with Adjustable Risks                                      
 }
%
\ARTICLEAUTHORS{%
\AUTHOR{Yiling Zhang}
\AFF{Industrial and Systems Engineering, University of Minnesota, \EMAIL{yiling@umn.edu}}
} 

	\ABSTRACT{%
	We study distributionally robust chance-constrained programs (DRCCPs)  with individual chance constraints under a Wasserstein ambiguity. The DRCCPs treat the risk tolerances associated with the distributionally robust chance constraints (DRCCs) as decision variables to trade off between the system cost and risk of violations by penalizing the risk tolerances in the objective function. The introduction of adjustable risks, unfortunately, leads to NP-hard optimization problems. We develop integer programming approaches for individual chance constraints with uncertainty either on the right-hand side or on the left-hand side. In particular, we derive mixed integer programming reformulations for the two types of uncertainty to determine the optimal risk tolerance for the chance constraint. Valid inequalities are derived to strengthen the formulations. We test diverse instances of diverse sizes.
	
	
}%


\HISTORY{}

\maketitle

%
\graphicspath{{figures/}}

\section{Introduction}
In many planning and operational problems, chance constraints are often used for ensuring the quality of service (QoS) or system reliability. For example, 
chance constraints can be used to restrict the risk of under-utilizing renewable energy in  power systems \citep[e.g.,][]{ma2019distributionally,zhang2022building}, to constrain the risk of loss in portfolio optimization \citep[e.g.,][]{lejeune2016multi}, and to impose the probability of satisfying demand in humanitarian relief networks \citep[e.g.,][]{elcci2018chance}.
%
In particular, with a predetermined risk tolerance $\alpha\in[0,1]$, a generic chance constraint is formulated in the following form.
\begin{equation*}\label{eq:cc}
	\mathbb P_{f}(T(\xi)x\ge q(\xi)) \ge 1-\alpha, 
\end{equation*}
where $x\in\mathbb R^d$ and the probability of violating the constraint $T(\xi)x\ge q(\xi)$ is no more than $\alpha$ with
a random vector $\xi\in\mathbb R^l$ following distribution distribution $f$. 
The technology matrix is obtained using a function $T(\xi):\mathbb R^l \mapsto\mathbb R^{m\times d}$ and the right-hand side (RHS) is a function $q(\xi):\mathbb R^l\mapsto \mathbb R^m$.

When an accurate estimate of the underlying distribution $f$ is not accessible, \emph{distributionally robust optimization} (DRO) provides tools to accommodate incomplete distributional information. Instead of assuming a known underlying distribution, DRO considers a prescribed set $\mathcal D$ of probability distributions, termed as an \emph{ambiguity set}. The distributionally robust variant of chance constraint \eqref{eq:cc} is as follows.
\begin{equation*}
	\inf_{f\in\mathcal D}\mathbb P_{f}(T(\xi)x\ge q(\xi)) \ge 1-\alpha. 
\end{equation*}
In the 
distributionally robust chance constraint (DRCC), $\alpha$ represents the worst-case probability of violating constraints $T(\xi)x\ge q(\xi)$ with respect to the ambiguity set $\mathcal D$.


In many system planning and operational problems,
a higher value of the probability $1-\alpha$ can lead to potentially better customer \textcolor{black}{satisfactions} and/or a lower probability of unfavorable events.
However, a too-large $1-\alpha$  may lead to problem infeasibility and requires additional resources and operational costs \citep[e.g.,][]{ma2019distributionally}.
%
%
%
To find a proper balance between the cost and reliability objectives, alternatively, in this paper, we consider DRCC problems with an adjustable risk, where
the risk tolerance $\alpha$ is treated as a  variable.

In particular, with a variable risk tolerance $\alpha$, we consider 
%
\begin{subequations}\label{eq:adjust-drcc}
	\begin{eqnarray}
		z_0: = \min_{x\in\mathcal X,\alpha} && c^\top x + g(\alpha)\\
		\mbox{s.t.} && \inf_{f\in\mathcal D}\mathbb P_{f}\left(T(\xi)x\ge q(\xi)\right) \ge 1-\alpha \label{eq:drcc}\\
		&& \alpha \in [0,\overline{\alpha}], \label{eq:drcc_alpha}
	\end{eqnarray}
\end{subequations}
%
The risk tolerance $\alpha$ is upper bounded by a parameter $\overline{\alpha}< 1$. 
The parameter $\bar \alpha$ is predetermined and can be viewed as the most risk of unfavorable events that the decision maker is willing to take. The objective trades off between the system cost $c^\top x$ with $c\in\mathbb R^d$ and the (penalty) cost of allowed violation risk $g(\alpha): \ [0,\bar{\alpha}]\rightarrow \mathbb R_0^+$. \textcolor{black}{The risk cost function $g(\alpha)$ is assumed monotonically increasing in $\alpha$. 
	While outside of the paper's scope, the choice of the function $g(\alpha)$ is complicated and problem-dependent. For instance, in vertiport planning, the risk level $1-\alpha$ can be interpreted as service adoption rate. With a uniform price per ride, $g(\alpha)$ is a linear function of the adoption rate.}

In the following, we focus on  individual chance constraints, i.e., $m=1$: (1) with RHS uncertainty -- a fixed technology matrix $T(\xi) = T\in\mathbb R^{d}$ and a random RHS $q(\xi) = \xi\in\mathbb R$ following a univariate distribution ($l=1$); and (2) with left-hand side (LHS) uncertainty -- a random technology matrix $T(\xi):\ \mathbb R^l \mapsto\mathbb R^{d}$ and a fixed RHS $q(\xi) = q \in\mathbb R$. For the LHS uncertainty, we specify the technology matrix $T(\xi)$ by assuming that $T(\xi)$ are affinely dependent on $\xi$, i.e., $T(\xi) = \xi^\top \bar A + \bar b^\top$. Under this uncertainty setting, we can reformulate constraint $T(\xi)x \ge q$ as follows:
$$T(\xi)x \ge q \Leftrightarrow \xi^\top \bar A x  +  \bar b^\top x\ge q \Leftrightarrow A(x)\xi \ge b(x),$$
where scalar function $b(x) = q - \bar b^\top x$ and $A(x) = x^\top \bar A^\top$ is a $1\times l$ vector.

In thise paper, we focus on the risk-adjustable formulation \eqref{eq:adjust-drcc} with a single individual chance constraint for clarity.	We note that the results in the paper are readily to extend to the variant with a number $\Theta$ of individual chance constraints (see, e.g., the computational study on the transportation problem in Section \ref{sec:comp}), i.e.,
\begin{subequations}\label{eq:adjust-drcc-multi}
	\begin{eqnarray}
		\min_{x\in\mathcal X,\alpha} && c^\top x + \sum_{\theta=1}^\Theta g_\theta (\alpha_\theta)\\
		\label{eq:adjust-drcc-multi-constr1}		\mbox{s.t.} && \inf_{f\in\mathcal D_\theta}\mathbb P_{f}\left(T_\theta (\xi)x\ge q_\theta (\xi)\right) \ge 1-\alpha_\theta, \ \theta =1,\ldots,\Theta \\
		\label{eq:adjust-drcc-multi-constr2}		&& (\alpha_1, \alpha_2, \ldots, \alpha_\Theta) \in W,
	\end{eqnarray}
\end{subequations}
which is NP-hard with a Wasserstein ambigutiy set (see Section \ref{sec:wass_np} for the definition of the Waserstein ambiguity set and a proof).
Here, $W \subset \mathbb [0,1]^\Theta$ is a polyhedral set.	For the RHS uncertainty case, $T_\theta (\xi) = T_\theta \in\mathbb R^{d}, \  q_\theta (\xi) = \xi_\theta \in\mathbb R, \ \theta =1,\ldots,\Theta $; for the LHS uncertainty case, $T_\theta (\xi) = \xi^\top \bar A_\theta + \bar b_\theta^\top : \ \mathbb R^l \mapsto\mathbb R^{d}$, $q_\theta (\xi) = q_\theta \in\mathbb R, \ \theta=1,\ldots,\Theta.$ 
\subsection{Motivation and Related Literature}
\label{sec:motivation}

The idea of using variable risk tolerances, {dating} back to \citet{evers1967new}, explores trade-off between costs and the probability of not meeting specifications in metal melting furnace operations. It later has been extensively applied across various fields, including facility sizing \citep{rengarajan2009estimating}, flexible ramping capacity \citep{wang2018adjustable}, power dispatch \citep{qiu2016data,ma2019distributionally}, portfolio optimization \citep{lejeune2016multi}, humanitarian relief network design \citep{elcci2018chance},  inventory control problem \citep{rengarajan2013convex}, and urban mobility planning \citep{kai2022vertiport}, to name a few. \textcolor{black}{In these domains, risk tolerance -- often interpreted as reliability, quality of service (QoS), or service adoption rate -- is treated as a decision variable to balance operational or system design cots against performance outcomes.}

\textcolor{black}{When the RHS is deterministic and uncertainty only resides on the LHS, optimizing DRCCs  with variable risk tolerance aligns with the spirit of the target-oriented decision making: rather than optimizing an objective with a predetermined risk tolerance, decision makers specify an acceptable target (i.e., the deterministic RHS) and minimize the risk of failing to meet the target. This approach,  known as \emph{target satisficing}, was introduced by \citet{simon1959decision}. When distributional ambiguity is present, \citet{long2023robust} propose a \emph{robust satisficing} model, extending the concept of satisficing to account for ambiguity in distributional information. The robust satisficing approach has later been applied in various problems, including $p$-hub center problems \citep{zhao2023distributionally}, supply chain performance  optimization \citep{chen2022supply}, bike sharing \citep{simon1959decision}, and optimal sizing in power systems \citep{keyvandarian2023optimal}. 
	Although developed independently from target satisficing, \citet{xu2012optimization} 	propose an extension by assuming the risk tolerance as a function of the target in probabilistic
	envelope constraints.}

\citet{rengarajan2009estimating, rengarajan2013convex} perform Pareto analyses to seek an efficient frontier for a trade-off between the total investment cost and the probability of disruptions that cause undesirable events. In particular, they  solve a series of chance-constrained programs for a large number of risk-level $\alpha$ choices.   The chance constraint is required to be met for a range of the target.
Unlike \citet{rengarajan2009estimating, rengarajan2013convex}, another stream of research treats the risk tolerance $\alpha$ as a decision variable and develops nonparametric approaches to trade off the cost and reliability. 
With only right-hand side uncertainty, \citet{shen2014using} develops a mixed integer linear programming (MILP) reformulation
for individual chance constraints with only RHS uncertainty ($m=1$ in the chance constraint \eqref{eq:cc}) under discrete distributions. Along the same line, \citet{elcci2018chance} propose an alternative MILP reformulation for the same setting using knapsack inequalities.
In the context of joint chance constraints ($m>1$), \citet{lejeune2016multi} use Boolean modeling framework to develop exact reformulations for the case with RHS uncertainty and inner approximations for the case with LHS uncertainty. 
All these studies assume known underlying (discrete) probability distributions. In recent work, \citet{zhang2022building}  consider the distributionally robust variants of the risk-adjustable chance constraints under ambiguity sets with moment constraints and Wasserstein metrics, respectively.  They consider an individual DRCC with RHS uncertainty. For the moment-based ambiguity set, they develop two second-order cone programs (SOCPs) with $\alpha$ in different ranges; for the Wasserstein ambiguity set, they propose an exact MILP reformulation. 
In particular,  the Wasserstein ambiguity set assumes a continuous support and their results
require decision variables to be all pure binary to facilitate the linearization of bilinear terms.  Without assuming pure binary decisions, in this paper, we will develop integer approaches for solving the DRCCs with adjustable risks under Wasserstein metrics.


\textcolor{black}{Another motivation for focusing on   individual chance constraints with adjustable risks is to provide an approximation of joint chance constraints $\inf_{f\in\mathbb D} \mathbb P_f (T_\theta(\xi)x \ge q_\theta(\xi), \ \theta=1,\ldots,\Theta) \ge 1-\alpha$. 
	This approximation specifies set $W=\left\{ \sum_{\theta=1}^\Theta \alpha_\theta \le \alpha\right\}$ in \eqref{eq:adjust-drcc-multi}  to the multiple chance constrained variant \eqref{eq:adjust-drcc-multi} and serves as a natural extension of the classic Bonferroni approximation, which fixed $\alpha_\theta$ with equal tolerance $\alpha/\Theta$.  Optimizing those risk tolerances $\alpha_\theta$ of individual chance constraints potentially leads to better approximations \citep[see, e.g.,][]{prekopa2003probabilistic,xie2022optimized}.}

\textcolor{black}{This work is also related to an increasing number of recent works on chance constraints under distributional ambiguity. Among them, \citet{chen2022data,xie2018distributionally} and \citet{ji2021data} discuss and provide tractable formulations for various chance constraints under Wasserstein ambiguity. As a follow-up, \citet{chen2023approximations} develop approximations based on conditional value-at-risk and Bonferroni's inequality.
	Although this paper will focus on Wasserstein ambiguity, we also acknowledge works concerning other types of ambiguity in chance constraints. For instance, \citet{jiang2013data} investigate tractable reformulations for DRCCs under a set of distributions based on a general $\phi$-divergence measure. \citet{xie2018deterministic} study deterministic reformulations for joint chance constraints under moment ambiguity. Under generalized moment bounds and structural properties, \citet{hanasusanto2015distributionally} present a unifying framework for solving DRCCs and uncertainty quantification problems. We refer interested readers to \citet{kuccukyavuz2022chance} for a recent review on chance constraints under distributional ambiguity.}


\subsection{Main Contributions}

In this paper, we study the risk-adjustable DRCC when the ambiguity set is specified as a Wasserstein ambiguity set.	
The main contributions of the paper are three-fold. 
First, we establish that optimizing the risk-adjustable DRCC with multiple chance constraints is strongly NP-hard. 
Second, we develop tractable integer programming approaches to solving the risk-adjustable DRCC with randm RHS and LHS, respectively:
(1) By exploiting the (hidden) discrete structures of the individual DRCC with random RHS, we provide tractable mixed-integer reformulations for risk-adjustable DRCCs using the Wasserstein ambiguity set.
%
Specifically, a MILP reformulation is proposed under the finite distribution assumption, and a mixed-integer second-order cone programming (MISOCP) reformulation is derived under the continuous distribution assumption.
Moreover, we strengthen the proposed mixed-integer reformulations by deriving valid inequalities by exploring the mixing set structure of the MILP reformulation and submodularity in the MISOCP reformulation. (2) With random LHS, we provide an equivalent mixed-inter conic reformulation when the decision $x$ is binary and a valid inequality is derived to improve the computation.
Third,  extensive numerical studies are conducted to demonstrate the computational efficacy of the proposed solution approaches.

The remainder of the paper is organized as follows. 
Section \ref{sec:wass_prelim} presents the Wasserstein ambiguity set and preliminary results regarding individual DRCC. Section \ref{sec:risk-adjustable_form} studies the case with RHS uncertainty and utilizes hidden discrete structures to derive mixed integer programming reformulations along with valid inequalities to strengthen the mixed integer programming reformulations under discrete and continuous distribution assumptions. Section \ref{sec:lhs} focuses on the case with LHS uncertainty and binary decisions. An equivalent mixed-inter conic reformulation is derived with a vali inequality.
Section \ref{sec:comp} demonstrates the computational efficacy of the proposed approaches for solving a transportation problem with diverse problem sizes and a demand response management problem. Finally, we draw conclusions in Section \ref{sec:conclusions}.

\section{Wasserstein Ambiguity and Preliminary Results}
\label{sec:wass_prelim}
\subsection{Wasserstein Ambiguity and NP-hardness}
\label{sec:wass_np}
In the risk-adjustable DRCC \eqref{eq:drcc}, we consider a Wasserstein ambiguity set $\mathcal D$ constructed as follows.
Given a series of  $N$ historical data samples $\{\xi^n\}_{n=1}^N$ drawn from  $\mathbb R^l$, the empirical distribution is constructed as $\mathbb P_0(\tilde \xi = \xi^n) = 1/N, \ n=1,\ldots,N$.
For a positive radius $\epsilon > 0$, the Wasserstein ambiguity set defines a ball around a reference distribution (e.g., the empirical distribution) in the space of probability distributions as follows:
\begin{equation*}
	\mathcal D: = \left\{f: \ \mathbb P_f(\tilde \xi \in \mathbb R^l)=1, \ W(\mathbb P_f,\mathbb P_0) \le \epsilon \right\}.
\end{equation*}
The Wasserstein distance is defined as
\begin{equation*}
	W(\mathbb P_f, \mathbb P_0) := \inf_{\mathbb Q\sim (\mathbb P_1,\mathbb P_2)}\mathbb E_{\mathbb Q} \left[\|\tilde \xi_1 - \tilde \xi_2 \|_p \right],
\end{equation*}
where $\tilde \xi_1$ and $\tilde \xi_2$ are random variables following distribution $\mathbb P_1$ and $\mathbb P_2$, $\mathbb Q\sim (\mathbb P_1, \mathbb P_2)$ denotes a joint distribution of $\tilde \xi_1$ and $\tilde \xi_2$ with marginals $\mathbb P_1$ and $\mathbb P_2$, and $\|\cdot\|_p$ denotes the $p$-norm. To exclude trivial special cases, throughout the rest of the paper, we assume that $\epsilon>0$ and \textcolor{black}{$\bar\alpha\in (0,1)$}. 

{\color{black}Under the Wasserstein ambiguity set,  solving the risk-adjustable DRCC problem is in general strongly NP-hard when multiple chance constraints are involved as in \eqref{eq:adjust-drcc-multi}.
	Denote $S = \left\{ (x,\alpha): \ \eqref{eq:adjust-drcc-multi-constr1}-\eqref{eq:adjust-drcc-multi-constr2}  \right\}$.
	
	\begin{theorem}\label{thm:np-hard}
		It is strongly NP-hard to optimize over set $S$.
		
	\end{theorem}
	
	\begin{proof}{Proof of Theorem \ref{thm:np-hard}:}
		
		We obtain this result by showing that the strongly NP-complete binary program can be reduced to verifying whether $(x,\alpha)$ belongs $S$.
		
		\emph{Binary programming feasibility}:	Given an integer matrix $P\in\mathbb Z^{\tau\times d}$ and integer vector $h\in \mathbb Z^\tau$, is there a solution $x\in\mathbb \{0,1\}^d$ such that  $Px \ge h$?
		
		To answer the binary programming feasibility problem given an instance with matrix $\hat P\in\mathbb Z^{\tau\times d}$ and vector $\hat h\in\mathbb Z^\tau$, we can verify whether the following instance of $S$ has a solution.
		\begin{equation}\label{eq:s_inst}
			S = \left\{(x,\alpha): \ \begin{array}{l}
				\inf_{f_i\in\mathcal D_i} \mathbb P_{f_i} \left(x_i \ge \xi_i \right) \ge 1-\alpha_i, \ i=1,\ldots,d\\
				\inf_{f_i\in\mathcal D_i} \mathbb P_{f_i} \left(1 - x_{i-d} \ge \xi_i \right) \ge 1-\alpha_i, \ i=d+1,\ldots,2d\\
				\inf_{f_i\in\mathcal D_i} \mathbb P_{f_i} \left(\hat P_{i-2d}x \ge \hat h_{i-2d} \right) \ge 1-\alpha_i, \ i=2d+1,\ldots,2d+\tau\\
				\sum_{i=1}^{2d+\tau} \alpha_i \le \frac{dK}{N}
			\end{array} \right\},\end{equation}
		where the Wasserstein ambiguity set $\mathcal D_i$ has a zero radius and contains a singleton, $\mathcal D_i = \{ f_0 \}$.
		The DRCCs in \eqref{eq:s_inst} reduce to stochastsic chance constraints with the full knowledge of the distribution  $f_0$.
		Here, $f_0$ is an empirical distribution constructed with samples drawn from a Bernoulli distribution. Assume that these samples are in a non-drecreasing order with the first $K<N$ samples of ones and the rest of zeros. 
		
		The first set of DRCCs in \eqref{eq:s_inst} is equivalent to
		\begin{equation}\label{eq:s_inst_drcc1}
			\alpha_i \ge 1 - \frac{1}{N}\left(K\mathbb{I} (x_i \ge 1)+(N-K)\mathbb{I} (x_i \ge 0)\right), \ i=1,\ldots,d,
		\end{equation}
		where $\mathbb I$ is the indicator function. The equivalence is due to that $\mathbb P_{f_0}(x_i\ge \xi_i) = \left(K\mathbb I(x_i\ge 1) + (N-K)\mathbb I(x_i\ge 0)\right)/N$. Similarly, the second set of DRCCs in \eqref{eq:s_inst} is equivalent to
		\begin{equation}\label{eq:s_inst_drcc2}
			\alpha_i \ge 1 - \frac{1}{N}\left(K\mathbb{I} (x_i \le 0)+(N-K)\mathbb{I} (x_i \le 1)\right), \ i=d+1,\ldots,2d.
		\end{equation}
		Following Fourier-Motzkin elinimiation of variables $\alpha_i, \ i=2d+1,\ldots,2d+\tau$, \eqref{eq:s_inst} is equivalent to
		\begin{equation}\label{eq:s_inst_fm}
			S = \left\{(x,\alpha): \ \eqref{eq:s_inst_drcc1}-\eqref{eq:s_inst_drcc2}, \ \hat Px \ge \hat h, \ \sum_{i=1}^{2d} \alpha_i \le \frac{dK}{N}  \right\}.
		\end{equation}
		A feasible $x_i$ must belong to  $[0,1]$. To this see, let us consider that (1) if $x_i > 1$, then following \eqref{eq:s_inst_drcc2}, we have  $\alpha_i\ge1, \ i=d+1,\ldots,2d$ and (2) if $x_i <0$,  \eqref{eq:s_inst_drcc1}  leads to $\alpha_i \ge 1, \ i=1,\ldots,d$. In both cases, constraint $\sum_{i=1}^{2d} \alpha_i \le dK/N$ in \eqref{eq:s_inst_fm} is violated and thus $x_i\in [0,1], \ i=1,\ldots, 2d$. Denote index sets $I_0 = \{i: \ x_i = 0, \ i=,1\ldots,d\}$, $I_1 = \{i: \ x_i = 1, \ i=,1\ldots,d\}$, and $I_2 = \{i: \ 0< x_i < 1, \ i=,1\ldots,d\}$. It is clear that $|I_0|+|I_1|+|I_2| = d$. We also have that
		\begin{subequations}
			\begin{eqnarray*}
				\sum_{i=1}^{2d}\alpha_i && \ge 2d - \frac{1}{N}\sum_{i=1}^d\left( K\mathbb{I} (x_i \ge 1)+(N-K)\mathbb{I} (x_i \ge 0) + K\mathbb{I} (x_i \le 0)+(N-K)\mathbb{I} (x_i \le 1) \right)\hspace{5mm}\\
				&& =2d - \frac{K|I_1| + 2(N-K)(|I_0|+|I_1|+|I_2|) + K|I_0|}{N}\\
				&& =2d - \frac{(2N-K)(|I_0|+|I_1|+|I_2|) - K|I_2|}{N} = \frac{K}{N}(d + |I_2|),
			\end{eqnarray*}
		\end{subequations} 
		where the first inequality is due to constraints \eqref{eq:s_inst_drcc1}-\eqref{eq:s_inst_drcc2} and the first equality holds following the definition of the three index sets. Given that $\sum_{i=1}^{2d} \alpha_i \le dK/N$, we have $|I_2| =0$ and thus $x\in\{0,1\}^d$. Consequently, $\alpha_i = (1-\mathbb I(x_i=1))K/N, \ i=1,\ldots,d$ and $\alpha_i = (1-K\mathbb I(x_i=0))K/N, \ i=d+1,\ldots,2d$. Now, we have that  the binary programming problem is feasible if and only if the instance \eqref{eq:s_inst} of $S$ has a feasible solution.

\end{proof}}

\subsection{Preliminary Results for DRCC with a Known Risk Tolerance $\alpha$}
\label{sec:prelim}

We denote $S(x) = \left\{\xi\in\mathbb R^l: \ T(\xi)x > q(\xi)\right\}$. Let $\text{int}S(x)$ denote its interior and $\text{cl}S(x)$ denote its closure.  Proposition 3 in \citet{gao2023distributionally} implies that $\inf_{f\in\mathcal D} \mathbb P(\xi \in \text{int}S(x)) = \inf_{f\in\mathcal D} \mathbb P(\xi \in S(x)) =\inf_{f\in\mathcal D} \mathbb P(\xi \in \text{cl}S(x)) = \inf_{f\in\mathcal D} \mathbb P(T(\xi)x\ge q(\xi))$. That is, regardless of whether $S(x)$ is open or closed, the worst-case probability remains unchanged.
With this, in what follows, we use the open set $S(x)$ for convenience and  denote $\bar S(x) = \mathbb R^l\backslash S(x)$ its closed complement.

In the rest of the paper, without loss of generality, we assume that the samples are ordered in non-decreasing distance to $\bar S$, i.e., $\mathrm{dist}(\xi^1,\bar S(x)) \le \mathrm{dist}(\xi^2, \bar S(x)) \le \cdots \le \mathrm{dist}(\xi^N, \bar S(x))$, where the distance  $\mathrm{dist}(\xi,\bar S(x)):= \min_{\xi\in \mathbb R^l} \left\{ \|\xi - \xi^\prime\|: \ \xi^\prime \in \bar S(x)\right\}$ with a a general norm $\|\cdot\|$. In the RHS uncertainty case, the distance is given by $\mathrm{dist}(\xi,\bar S(x)) = (Tx -\xi)^+$, where operator $(a)^+ = \max \{a,0\}$ takes the positive part of $a$. In the LHS uncertainty case, the distance is given by  $\mathrm{dist}(\xi,\bar S(x)) =(A(x)\xi^{(n)} - b(x))^+/\|A(x)\|_*$ (e.g., Lemma A.1 in \citet{chen2022data}), where $\|\cdot\|_*$ represents the dual normal of $\|\cdot\|$.

\begin{proposition}[Adapted from Theorem 2 in \citet{chen2022data}]
\label{prop:dist}
For a given risk tolerance $\alpha$, the DRCC \eqref{eq:drcc} is equivalent to 

\begin{equation}
	\sum_{n=1}^{\alpha N} \mathrm{dist}(\xi^{(n)},\bar S(x)) \ge N\epsilon, \text{ or } \max\left\{j\in [0,N] \bigg|\sum_{n=1}^{j} \mathrm{dist}(\xi^{(n)}, \bar S(x))\le \epsilon \right\}\le \alpha N,\label{eq:drcc-dist}
\end{equation}
where the summation in the first constraint is a partial sum for fractional $\alpha N$: $\sum_{n=1}^{\alpha N} k_n = \sum_{n=1}^{\lfloor \alpha N \rfloor}k_n + (\alpha N - \lfloor \alpha N \rfloor) k_{\lfloor\alpha N \rfloor+1}$. 

\end{proposition}

\begin{figure}[htb]
\centering
\includegraphics[width=0.48\linewidth]{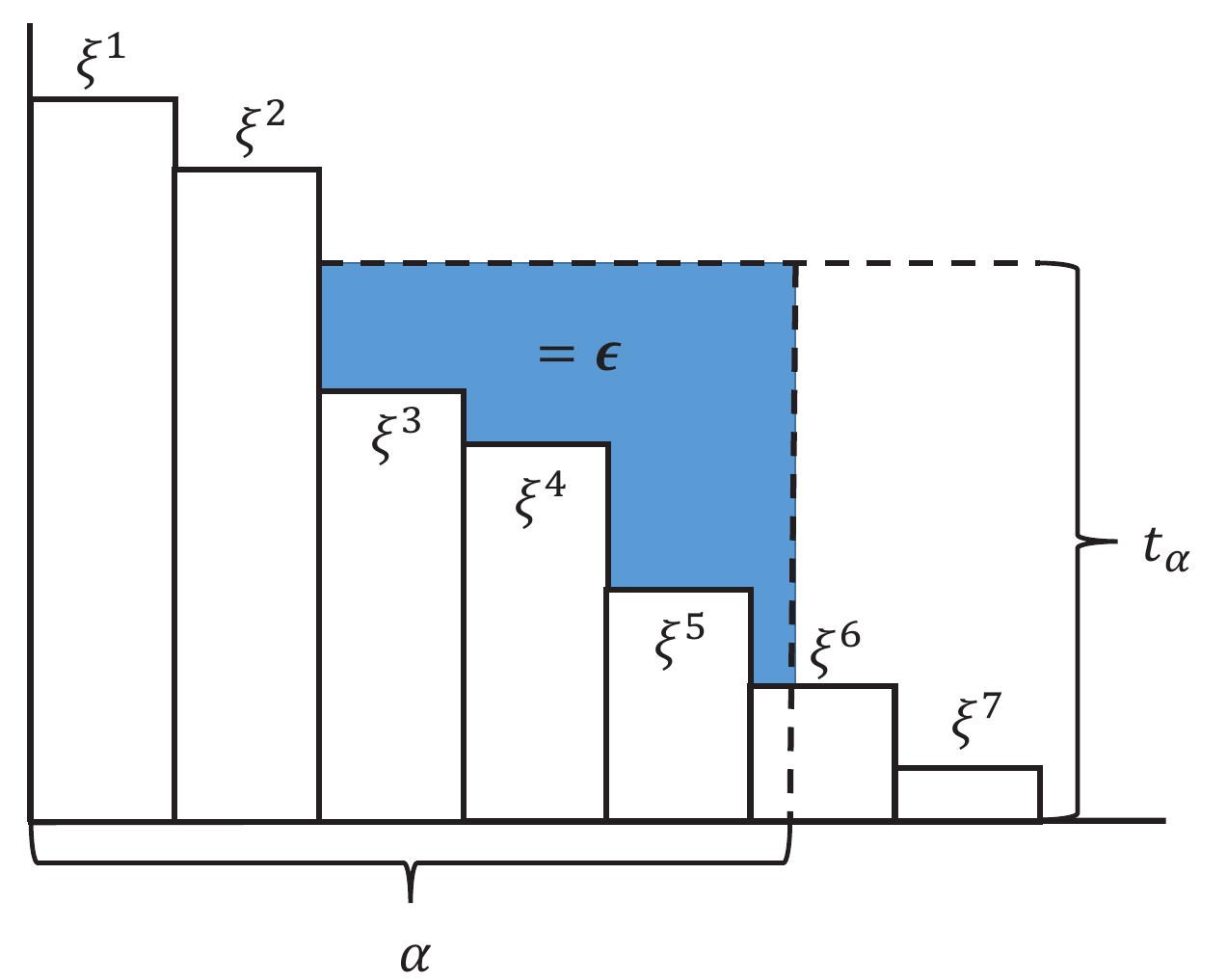}
\caption{Illustration of the water-filling interpretation for the partial-sum inequality in \eqref{eq:drcc-dist} with RHS uncertainty.}
\label{fig:water-filling}
\end{figure}

With only RHS uncertainty, the non-decreasing distance order implies a non-increasing order of samples, i.e., $\xi^1 \ge \xi^2\ge \cdots \ge \xi^N$. .
The first equivalent constraint in \eqref{eq:drcc-dist} has a \emph{water-filling} interpretation as illustrated in Figure \ref{fig:water-filling}. The height of patch $n$ is given by $\xi^n$ and the width is given by $1/N$. The region with a width of $\alpha$ is flooded to a level $t_\alpha$ which uses a total amount of water equal to $\epsilon$. Then the reformulation of DR chance constraint \eqref{eq:drcc-dist} is equivalent to a linear inequality $Tx\ge t_\alpha$. The water level $t_\alpha$ represents the worst-case value-at-risk (VaR):
\begin{equation}\label{eq:t-def}t_\alpha := \inf_v\left\{v: \ \inf_{f\in\mathcal D}\mathbb P(v\ge \xi) \ge 1-\alpha \right\}  = \min_v\left\{ v: \ \frac{1}{N}\sum_{n=1}^{\alpha N} (v-\xi^n)^+\ge \epsilon \right\}.\end{equation}

Let $j^*$ be the largest index such that when  the amount of water that fills the region of width $\alpha$ to the level $\xi^{j^*}$ is no less than $\epsilon$. That is, \begin{equation}\label{eq:dist_j}j^*: = \max\left\{j\in\{1,\ldots,N\}: \ \xi^j \ge \frac{N\epsilon + \sum_{n=j+1}^{\alpha N}\xi^n }{N\alpha - j} \right\}.\end{equation}
For example, in Figure \ref{fig:water-filling}, if the water is filled up to the level as $\xi^1$ or $\xi^2$, the amount of water exceeds $\epsilon$. In this example, $j^* = 2$.
We note that such index $j^*$ may not always exist, i.e., the  problem  \eqref{eq:dist_j} can be infeasible. This happens if the amount of water is strictly less than $\epsilon$ even when the water level reaches $\xi^1$.
In this case, one can keep increasing the water level until the amount of water equals $\epsilon$ and let the worst-case VaR $t_\alpha$ equal the water level. Otherwise, the worst-case VaR can be obtained using the  propositions below.

\begin{proposition}[Adapted from Theorem 2 in \citet{ji2021data}]
When the random RHS $\xi$ has a finite support with unknown mass probability, the worst-case VaR $t^d_\alpha = \xi^{j^*}.$  
\label{prop:finite_dis_var}
\end{proposition}

\begin{proposition}[Adapted from Theorem 8 in \citet{ji2021data}]
When  the random RHS $\xi$ has a continuum of realizations, the worst-case VaR $$t^c_\alpha = \frac{N\epsilon + \sum_{n=j^*+1}^{\alpha N}\xi^n }{N\alpha - j^*}.$$
\label{prop:conti_dis_var}
\end{proposition}
It is easy to verify that $t_\alpha^d = \xi^{j^*} \ge t_\alpha^c$ given the definition of the critical index $j^*$ in \eqref{eq:dist_j}.


%
%
%
%
%
%

\section{Risk-Adjustable DRCC with Right-Hand Side Uncertainty}
\label{sec:risk-adjustable_form}
In this section, 
we consider two uncertainty types for the risk-adjustable DRCC model \eqref{eq:adjust-drcc} with RHS uncertainty ($T(\xi)=T, \ q(\xi) = \xi$): 
\begin{enumerate}[label={A\arabic*}]

\item \emph{Finite Distribution}: The random vector $\xi$ has a finite support. The mass probability of each atom is unknown and allowed to vary.
\label{assump:finite}

\item \textcolor{black}{\emph{Continuous Distribution}: The random variable $\xi$ has a continuum (infinite number) of realizations and the probability of every single realization is zero.} \label{assump:continuous}
\end{enumerate}
%
%
%
First, in Section \ref{sec:relation_finite_continuous}, we provide the relation of the optimal values under the two distribution assumptions. Then, we 
develop tractable mixed-integer reformulations and valid inequalities under the finite and continuous distributions in Sections \ref{sec:finite_dist} and \ref{sec:continuous}, respectively. 
In Appendix \ref{sec:dominance}, we generalized the concept of non-dominated points, or the so-called $p$-efficient points to the distributionally robust setting.

\subsection{Relation of the Optimal Values under the Two Distribution Assumptions}
\label{sec:relation_finite_continuous}

Consider a function $t(\alpha)$ which maps the risk tolerance $\alpha$ to its corresponding worst-case VaR. If the function is known, the DRCC \eqref{eq:drcc} is equivalent to a linear constraint of $x$.
Thus, the risk-adjustable DRCC problem \eqref{eq:adjust-drcc} is rewritten as follows.
\begin{equation}\label{eq:linear_equ}
z(t(\alpha)):=\min_{x\in\mathcal X, \alpha\in [0,\bar\alpha]}\left\{ c^\top x+ g(\alpha): \ Tx\ge t(\alpha) \right\},
\end{equation}
where the optimal value depends on the choice of function $t(\alpha)$. Under Assumption \ref{assump:finite} of the finite distribution, let $t^d(\alpha)$ be the worst-case VaR function and  the optimal value of \eqref{eq:linear_equ} be $z^d: = z(t^d(\alpha))$. Similarly, under Assumption \ref{assump:continuous} of the continuum realizations, let $t^c(\alpha)$ be the worst-case VaR function and the optimal value of \eqref{eq:linear_equ} be $z^c: = z(t^c(\alpha))$. 
The next proposition presents the relation between the optimal values with finite and continuous distributions.
%
\begin{proposition}\label{prop:obj}
The risk-adjustable DRCC problem \eqref{eq:adjust-drcc} under the continuous distribution assumption \ref{assump:continuous} yields an optimal value no more than that under the finite distribution assumption \ref{assump:finite}, i.e., $z^d\ge z^c$.
\end{proposition}
The proposition is an immediate result from the fact that, for a given $\alpha$, $t_\alpha^d\ge t_\alpha^c$.

\subsection{Finite Distribution}
\label{sec:finite_dist}


According to Propositions \ref{prop:finite_dis_var} and \ref{prop:conti_dis_var}, the worst-case VaR $t_\alpha^d$ (if exists) under the finite distribution assumption is the smallest $\xi^j$ which is no less than the worst-case VaR $t_\alpha^c$ under the continuous distribution assumption. That is, $t_\alpha^d = \min_{j\in \{1,\ldots,N\}} \left\{\xi^j:\  \xi^j\ge t_\alpha^c \right\}$.
We thus have the following result, which has already been anticipated in Proposition \ref{prop:finite_dis_var}.
\begin{corollary}\label{coro:finite_dist_lp}
Under the finite distribution assumption \ref{assump:finite} , for any risk tolerance $\alpha\in(0,1)$ such that $\xi^j \ge t_\alpha^c> \xi^{j+1} $ for some $j\in \{1,\ldots,N\}$, the DRCC $ \inf_{f\in\mathcal D}\mathbb P_{f}(Tx\ge \xi) \ge 1-\alpha$ is equivalent to a linear constraint:
$$Tx\ge \xi^j.$$
\end{corollary}
Given any fixed $\alpha_1,\alpha_2\in (0,1)$ such that $\xi^j \ge t_{\alpha_1}^c\ge t_{\alpha_2}^c > \xi^{j+1}$, the DRCCs of the two risk tolerances yield the same linear reformulation $Tx\ge \xi^j$ under the finite distribution assumption. Thus, in the risk-adjustable DRCC problem \eqref{eq:adjust-drcc}, it suffices to strengthen
$\alpha\in(0,\bar\alpha]$ by restricting it to the risk tolerances $\alpha$ such that
the corresponding worst-case VaR $t_\alpha^d\in(0,1)$ belongs to a discrete set:
$$t_\alpha^d\in \{\xi^1,\ldots, \xi^N\}.$$
\textcolor{black}{With this observation, when only a single DRCC is involved, one may solve at most $N$ problems, each replacing the DRCC with a linear constraint with one realization of $\xi$ from the discrete set, and select the one with minimum objective. However, when dealing with multiple DRCCs, the number of problems one needs to solve grows exponentially in $N$. Next, we develop an integer programming formulation, of which the size grows linearly in $N$.}

For each sample $\xi^n, \ n=1,\ldots,N$ in the discrete set, a risk tolerance $\alpha_n$, which achieves the worst-case VaR at $t^d_{\alpha_n} = \xi^n$, can be obtained using a bisection search method as the flooded area is non-decreasing in the risk tolerance (see the water-filling interpretation in Section \ref{eq:adjust-drcc}). 
We note that when $\xi^n$ is too small, the corresponding risk tolerance may not exist.
In this section, we assume that  the corresponding risk tolerances exist for the first $N^\prime\le N$ largest  samples $\left\{ \xi^n \right\}_{n=1}^{N^\prime}$ and denote their corresponding  risk tolerances by $\alpha_n, \ n=1,\ldots,N^\prime$. 
\begin{theorem}\label{thm:mip_discrete}
Under the finite distribution assumption, the risk-adjustable DRCC problem \eqref{eq:adjust-drcc} is equivalent to the following MILP formulation.
\begin{subequations}\label{eq:mip_discrete}
\begin{eqnarray}
\label{eq:mip_discrete_obj} z^d = \min_{x\in\mathcal X, y} && c^\top x + \sum_{n=1}^{N^\prime-1}\Delta_n y_n + \Delta_{N^\prime}  \\
\label{eq:mip_discrete_constr1}\mbox{s.t.} && Tx \ge \xi^n - M_n(1-y_n), \ n=1,\ldots,N^\prime-1\\
\label{eq:mip_discrete_alpha}&& \sum_{n=1}^{N^\prime-1}(\alpha_{n}-\alpha_{n+1})y_{n} + \alpha_{N^\prime} \in (0,\bar{\alpha}]\\
\label{eq:mip_discrete_binary}&& y_n\in\{0,1\}, \ n=1,\ldots,N^\prime-1,
\end{eqnarray}
\end{subequations}
where $\Delta_{n}:=g(\alpha_{n}) - g(\alpha_{n+1})\le 0, \ n=1,\ldots,N^\prime-1$, $\Delta_{N^\prime} := g(\alpha_{N^\prime})$ and $M_n$ is a big-M constant. 
\end{theorem}
\begin{proof}{Proof of Theorem \ref{thm:mip_discrete}:}
To see the equivalence, we need to show (1) $z^d\le z_0$ and (2) $z^d \ge z_0$. Recall that $z_0$ is the optimal value of the risk-adjustable DRCC problem \eqref{eq:adjust-drcc}.
\begin{enumerate}
\item[(1)] $z^d\le z_0$: Given an optimal solution $( x_0, \alpha_0)$ to the risk-adjustable DRCC problem \eqref{eq:adjust-drcc}, we will construct a feasible solution to MILP \eqref{eq:mip_discrete}.  Let $\bar y_n = 1$ if $T x_0 \ge \xi^n$ and $\bar y_n=0$ otherwise, for $ n=1,\ldots,N^\prime$. Then the solution $(x_0,\bar y_{n}, n=1,\ldots,N^\prime)$ satisfy constraints \eqref{eq:mip_discrete_constr1} and \eqref{eq:mip_discrete_binary}.

Let $j^*$ be the smallest index such that $T x_0 \ge \xi^{j^*}$. We will show that $\alpha_0 = \alpha_{j^*}$.
When $j^* = 1$, $t_{\alpha_0}^d = \xi^1$ and $\alpha_0 = \alpha_1$. When $j^* \ge 2$, we prove by contradiction by assuming two cases (i) $t_{\alpha_0}^d > \xi^{j^*}$ and (ii) $t_{\alpha_0}^d < \xi^{j^*}$. In the first case, $t_{\alpha_0}^d \le \xi^{j^*-1}$. According to Proposition \ref{prop:VaR}, $\alpha_0 > \alpha_{j^*-1}$. Then, $(x_0,\alpha_{j^*-1})$ is feasible to the risk-adjustable DRCC \eqref{eq:adjust-drcc} with a smaller objective value than $z_0$ as the function $g(\alpha)$ is increasing in $\alpha$. In the second case, $\alpha_0 \ge \alpha_{j^*}$ due to Proposition \ref{prop:VaR} and $(x_0,\alpha_{j^*})$ is a feasible solution with a smaller objective than $z_0$ in the risk-adjustable DRCC \eqref{eq:adjust-drcc}. Both cases result in a contradiction to the fact that $z_0$ is the optimal value of the risk-adjustable DRCC \eqref{eq:adjust-drcc}.



Since $\alpha_0  = \alpha_{j^*}$ and $\alpha_0 \in (0,\bar\alpha]$,
$\alpha_{j^*} = \sum_{n=2}^{N^\prime}(\alpha_{n-1}-\alpha_n)\bar{y}_{n-1} + \bar{y}_{N^\prime}\alpha_{N^\prime}\in (0,\bar\alpha]$ satisfies constraint \eqref{eq:mip_discrete_alpha}. Solution $(x_0,\bar y)$ is feasible to \eqref{eq:mip_discrete} with $c^\top x_0 + g(\alpha_{j^*}) = z_0$. Thus, $z^d \le z_0$.
\item[(2)] $z^d\ge z_0$: Given an optimal solution $(\hat x,\hat y)$ to problem \eqref{eq:mip_discrete}, we construct a feasible solution to the risk-adjustable DRCC problem \eqref{eq:adjust-drcc}. Denote $j$ the smallest index such that $T\hat x \ge \xi^j$, or, equivalently, the smallest index such that $\hat y_j = 1$.
Let $\hat\alpha = \alpha_j$. It is easy to see that $(\hat x,\hat \alpha)$ is feasible to the risk-adjustable DRCC problem \eqref{eq:adjust-drcc} and its objective value is $c^\top \hat x + g(\alpha_j) = z^d$. So $z^d \ge z_0$.
\end{enumerate}
Combining the two statements above completes the proof.
\end{proof}
\begin{remark}
A tight bound on the big-M constant $M_n$ is  $\xi^n - \xi^{N^\prime}$.
\end{remark}
Recall that the non-increasing order of samples: $\xi^1\ge\xi^2\ge\cdots\ge \xi^N$.
Thus, given an optimal solution $\bar x$ to MILP \eqref{eq:mip_discrete}, there exists a threshold index $j^*$ such that $T\bar x \ge \xi^{i}$, for any $i\ge j^*$, and $T\bar x < \xi^i$, for any $i<j^*$. 
As the objective coefficient  $\Delta_n, \ n=1,\ldots,N^\prime -1$ in \eqref{eq:mip_discrete_obj} are non-positive, in the optimal solution $(\bar x,\bar y)$, we have $\bar y_n = 1$, for $i\ge j^*$, and $\bar y_n=0$ for $i>j^*$.
By exploiting this solution structure, the next proposition presents how 
the MILP formulation \eqref{eq:mip_discrete} can be strengthened.
\begin{proposition}\label{prop:finite_strengthening}
~
\begin{itemize}
\item[i.] 
The following inequalities are valid for the MILP \eqref{eq:mip_discrete}:        
\begin{equation}\label{eq:mip_discrete_y_order}
y_{n+1} \ge y_{n}, \ n=1,\ldots,N^\prime -1.\end{equation}
\item[ii.] The \emph{strengthened star inequality} \citep{luedtke2010integer} is valid for the MILP \eqref{eq:mip_discrete}:
\begin{equation}
\label{eq:mip_discrete_constr_star}Tx\ge \xi^{N^\prime} + \sum_{n=1}^{N^\prime - 1}(\xi^n - \xi^{n+1})y_n.
\end{equation}
\end{itemize}
\end{proposition}
\begin{proof}{Proof of Proposition \ref{prop:finite_strengthening}:}
The valid inequalities \eqref{eq:mip_discrete_y_order} follow from the discussion above.
To see the second statement of the extended star inequalities, we introduce binary variable $z_n = 1-y_n, \ n=1,\ldots,N^\prime-1.$ Without loss of generality, we assume that $\xi^n\ge0$. Constraints \eqref{eq:mip_discrete_constr1} and \eqref{eq:mip_discrete_alpha} lead us to consider a \emph{mixing set} \citep{atamturk2000mixed,gunluk2001mixing, luedtke2010integer}:
\begin{equation}
P = \left\{ (t,z)\in\mathbb R_+ \times \{0,1\}^{N^\prime - 1}: \ \sum_{n=1}^{N^\prime-1}(\alpha_{n}-\alpha_{n+1})y_{n} + \alpha_{N^\prime}\le \bar\alpha,  \ t+ z_n \xi^n \ge \xi^n, \ n=1,\ldots,N^\prime \right\}
\end{equation}
where $t = Tx$. According to Theorem 2 in \citet{luedtke2010integer}, constraint  \eqref{eq:mip_discrete_constr_star} is face-defining for $\text{conv(P)}$. The proof is complete.
\end{proof}
\begin{remark}
When the distribution is known, a similar formulation for the stochastic chance-constrained problem can also be derived based on the Sample Average Approximation \citep{luedtke2008sample}. In this case, let $\alpha_n$ be the allowed risk tolerance when the VaR equals $\xi^n$ and $\alpha_n = n/N$.
The detailed MILP formulation for the stochastic chance-constrained formulation can be found in Appendix  \ref{apx:cc}. We note that the MILP formulation in Appendix  \ref{apx:cc} can be viewed as a hybrid of those in \citet{shen2014using,elcci2018chance}.
\end{remark}



\subsection{Continuous Distribution}
\label{sec:continuous}

Unlike the case with finite distributions, under the continuous distribution assumption \ref{assump:continuous}, the worst-case VaR cannot be restricted to a discrete set.

For a given risk tolerance $\alpha$, 
constraint \eqref{eq:drcc-dist} is equivalent to
\begin{equation}\label{eq:soc}
(\alpha N-j)(Tx-\xi^{k+1}) - \sum_{i=j+1}^k(\xi^i-\xi^{k+1}) \ge N\epsilon 
\end{equation}
where $k = \lfloor \alpha N \rfloor$ and
$j$ is the smallest index such that $Tx - \xi^{j+1} \ge 0$. For instance, in Figure \ref{fig:water-filling}, $j=2$ and $k=5$. When the risk tolerance $\alpha$ is not known, we introduce a binary variable $o_{jk}\in \{0,1\}$ to indicate if $j$ and $k$ are the two critical indices. Denote $\xi^0$ be an upper bound of $\xi$. 
We consider a mild assumption: 
\begin{enumerate}[label={A3}] 
\item For an optimal solution $\hat x$, $T\hat x\ge \xi^N$. That is, the optimal solution is restricted by the smallest realization of $\xi$.  \label{assump:x_stric}
\end{enumerate}

\begin{theorem}\label{thm:mip_continuous}
Under the continuous distribution assumption \ref{assump:continuous} and Assumption \ref{assump:x_stric}, the risk-adjustable DRCC problem is equivalent to the following mixed 0-1 conic formulation.
\begin{subequations}\label{eq:mip-continuous}
\begin{eqnarray}
z^c=  \min_{x\in\mathcal X, o,\alpha,u,w}&& c^\top x + g(\alpha)\\
\label{eq:mip-continuous_conic}\mbox{s.t.} && uw \ge \sum_{j=0}^{N-1}\sum_{k=j}^{N-1}o_{jk}\sum_{i=j+1}^k(\xi^i-\xi^{k+1}) + N\epsilon \\
\label{eq:mip-continuous_conic1}&& u\le \alpha N-\sum_{j=0}^{N-1}\sum_{k=j}^{N-1}jo_{jk}\\
\label{eq:mip-continuous_conic2}            && w \le Tx - \sum_{j=0}^{N-1}\sum_{k=j}^{N-1} \xi^{k+1}o_{jk}\\
\label{eq:mip-continuous_Tx}      && \sum_{j=0}^{N-1}\sum_{k=j}^{N-1} \xi^{j}o_{jk} \ge Tx \ge \sum_{j=0}^{N-1}\sum_{k=j}^{N-1} \xi^{j+1}o_{jk}\\
\label{eq:mip-continuous_alphaN}      && \sum_{j=0}^{N-1} \sum_{k=j}^{N-1} (k+1) o_{jk} \ge \alpha N \ge \sum_{j=0}^{N-1} \sum_{k=j}^{N-1} k o_{jk}\\
\label{eq:mip-continuous_sos1}&& \sum_{j=0}^{N-1}\sum_{k=j}^{N-1} o_{jk}=1\\
&& \alpha \in (0,\bar{\alpha}]\\
&& w\ge 0, \ u\ge 0\\
&& o_{jk}\in \{0,1\}, \ 0\le j\le k \le N-1.
\end{eqnarray}
\end{subequations}
\end{theorem}
\begin{proof}{Proof of Theorem \ref{thm:mip_continuous}:}
To establish the equivalence, we first show that $z^c\le z_0$ by constructing a feasible solution to problem \eqref{eq:mip-continuous} given an optimal solution to the risk-adjustable DRCC problem \eqref{eq:adjust-drcc}. Let $(x_0,\alpha_0)$ be an optimal solution to \eqref{eq:adjust-drcc}. Denote $k^* = \lfloor \alpha_0 N \rfloor$ and $j^*$ as the smallest index such that $Tx_0 - \xi^{j^*+1}\ge 0$. Let $\bar o_{j^*k^*} = 1$, $\bar o_{jk}=0, \ j\neq j^*, \ k\neq k^*, \ 0\le j\le k\le N-1$, $\bar u = \alpha_0N+1-j^*$, and $\bar w=Tx_0 - \xi^{k^*+1}$. 
It is easy to verify that $(x_0,\bar o, \alpha_0,\bar u,\bar w)$ is a feasible solution and its objective value equals $z_0$. Thus, $z^c\le z_0$.

To see the opposite direction
$z^c \ge z_0$, consider an optimal solution $(\hat x,\hat o,\hat \alpha,\hat u,\hat w)$ to problem \eqref{eq:mip-continuous}. Since $\hat o$ is feasible, there exists $\hat o_{\hat{j}\hat{k}}$ such that $\hat o_{\hat{j}\hat{k}}=1$ and $\hat o_{jk}=0, \ j\neq\hat j,\ k\neq \hat k$.
Combining constraints \eqref{eq:mip-continuous_conic}--\eqref{eq:mip-continuous_conic2}, we obtain 
\begin{equation}\label{eq:mip-continuous_conic_eq}
(\hat\alpha N+1-\hat j)(T\hat x-\xi^{\hat k+1}) - \sum_{i=\hat j+1}^{\hat k}(\xi^i-\xi^{\hat k+1})\ge N\epsilon.\end{equation} 
Constraints \eqref{eq:mip-continuous_Tx} and \eqref{eq:mip-continuous_alphaN} are equivalent to  
\begin{equation}\label{eq:mip-continuous_Tx_alphaN}
\xi^{\hat j }\ge T\hat x \ge \xi^{\hat j +1} \text{ and } \hat k+1\ge \alpha N \ge \hat k,    
\end{equation}
respectively.
Constraints \eqref{eq:mip-continuous_conic_eq} and \eqref{eq:mip-continuous_Tx_alphaN} imply that $(\hat x,\hat \alpha)$ satisfies the DR chance constraint \eqref{eq:drcc}. Thus, $(\hat x,\hat \alpha)$ is feasible to the risk-adjustable DRCC problem \eqref{eq:adjust-drcc} and $z^c\ge z_0$ as expected.
\end{proof}
\begin{remark}
The mixed 0-1 conic reformulation \eqref{eq:mip-continuous} consists of $(N^2-N)/2$ (additional) binary variables and two continuous variables.
When the decision $x\in\mathcal X\subset \{0,1\}^d$ is restricted to binary variables, under the continuous distribution assumption, \citet{zhang2022building} propose a MILP formulation (details are in Appendix \ref{apx:milp}) by linearizing bilinear terms in the quadratic constraint \eqref{eq:soc} using McCormick inequalities \citep[see, e.g.,][]{mccormick1976computability}. 
In addition to $(N^2-N)/2$  binary variables as those in the conic reformulation \eqref{eq:mip-continuous},  the linearization  introduces $(N^2-N)(2d+1)$ continuous variables, where $d$ is the dimension of $x$. The MILP reformulation usually does not scale well when the problem size grows, partly due to the weaker relaxations caused by the big-M type constraints, and also due to a larger number of added variables and constraints. We will later show the computational comparison in Section \ref{sec:comp_continuous}.

\end{remark}

In the mixed 0-1 conic reformulation \eqref{eq:mip-continuous}, 	
\textcolor{black}{constraint \eqref{eq:mip-continuous_conic} is a rotated conic quadratic mixed 0-1 constraint due to the fact that $o_{jk} = o_{jk}^2$ for binary $o_{jk}$.}
%
%
Although the resulting mixed-integer conic reformulation can be directly solved by optimization solvers, mixed 0-1 conic programs are often time-consuming to solve, mainly due to the binary restrictions.
%
In the following, we will develop valid inequalities for the mixed 0-1 conic reformulation \eqref{eq:mip-continuous} to help accelerate the branch-and-cut algorithm for solving \eqref{eq:mip-continuous}. Specifically, we explore the \emph{submodularity} structure of constraint \eqref{eq:mip-continuous_conic} as follows.

We first note that constraint \eqref{eq:mip-continuous_conic} can be rewritten in the following form 
\begin{equation}
\label{eq:mip-continuous_conic_general}
\sigma + \sum_{j=0}^{N-1}\sum_{k=j}^{N-1}d_{jk}o_{jk} \le uw,
\end{equation}
where $\sigma = N\epsilon>0$ and $d_{jk} = \sum_{i=j}^k(\xi^i-\xi^{k+1})\ge 0, \ j=0,\ldots,N-1, \ k=j,\ldots,N-1$.
By introducing auxiliary variable $\tau\ge 0$, constraint \eqref{eq:mip-continuous_conic_general} is equivalent to 
\begin{subequations}
\label{eq:mip-continuous_conic_general_soc}
\begin{eqnarray}
\label{eq:mip-continuous_conic_general_soc1}
&&\sqrt{\sigma + \sum_{j=0}^{N-1}\sum_{k=j}^{N-1}d_{jk}o_{jk}} \le \tau\\
\label{eq:mip-continuous_conic_general_soc2}  && \sqrt{\tau^2 + (w-u)^2} \le w+u.
\end{eqnarray}
\end{subequations}
The two inequalities \eqref{eq:mip-continuous_conic_general_soc1}--\eqref{eq:mip-continuous_conic_general_soc2} above are two second-order conic (SOC) constraints.   In particular, 
the convex hull of the first constraint \eqref{eq:mip-continuous_conic_general_soc1} can be fully described utilizing extended polymatroid inequalities as the left-hand side of constraint \eqref{eq:mip-continuous_conic_general_soc1} is a \emph{submodular} function \citep[see, e.g.,][]{atamturk2008polymatroids,atamturk2020submodularity}.
\begin{definition}[Submodular Function] Define the collection of set $[(N^2-N)/2]$'s subsets $\mathcal C:= \{S: \ \forall S\subset [(N^2-N)/2] \}$. 
Given a set function $g$: $\mathcal C\rightarrow \mathbb R$, g is submodular if and only if
\begin{equation*}
g\left(S \cup  \{j\}\right) - g(S) \ge g\left(R \cup \{j\} \right) - g\left( R\right),
\end{equation*}
for all subsets $S\subset R \subset \mathcal C$ and all elements $j\in \mathcal C\backslash R.$
\end{definition}
We use $g(S)$ and $g(o)$ interchangeably, where $o\in\{0,1 \}^{(N^2-N)/2}$ denotes the indicating vector of $S\subset\mathcal C$, i.e., $o_s=1$ if $s\in S$ and $o_s=0$ otherwise. The left-hand side of constraint \eqref{eq:mip-continuous_conic_general_soc1}, $h(o) :=\sqrt{\sigma + \sum_{j=0}^{N-1}\sum_{k=j}^{N-1}d_{jk}o_{jk}} $ is a submodular function, where $o$ is a one dimensional vector consisting of $o_{jk}, \ 0\le j\le k\le N-1$. 

\begin{definition}[Extended Polymatroid]
For a submodular function $g(S)$, the polyhedron
$$EP_g = \left\{\pi\in\mathbb R^{(N^2-N)/2}: \ \pi(S)\le g(S), \ \forall S\subset \mathcal C \right\}$$
is called an extended polymatroid associated with $g$, where $\pi(S) = \sum_{i\in S}\pi_i$.
\end{definition}

For submodular function $h$, 
linear inequality 
\begin{equation}
\pi^\top o \le z\label{eq:extended_polymatroid_ineq}
\end{equation}
is valid for the convex hull of the epigraph of $h$, i.e., conv$\{(o,z)\in\{0,1\}^{(N^2-N)/2}\times\mathbb R: \ z\ge h(o) \}$, if and only if $\pi$ is in the extended polymatroid, i.e., $\pi\in EP_h$ \citep[see][]{atamturk2008polymatroids}. The inequality \eqref{eq:extended_polymatroid_ineq} is called \emph{extended polymatroid inequality}. 

Although it suffices to only impose the extended polymatroid inequality at the extreme points of the extended polymatroid $EP_h$, there are an exponential number of them. Instead of adding all of them to the formulation  \eqref{eq:mip-continuous}, one can add them as needed in a branch-and-cut algorithm. Moreover,  the separation of the valid inequality \eqref{eq:extended_polymatroid_ineq} can be done efficiently using a $O(n\log{n})$ time greedy algorithm as follows.
Given a solution $(\hat o,\hat z)\in [0,1]^{(N^2-N)/2}\times\mathbb R_+$, one can obtain a permutation $\{(1),\ldots,(N^2)\}$ such that the elements of $o$ are sorted in a non-increasing order, $o_{(1)}\ge \ldots \ge o_{(N^2-N)/2}$. Let $S_{(i)} := \{(1),\ldots, (i) \}, \ i=1,\ldots, (N^2-N)/2$. Calculate $\hat\pi_{(1)} = h(S_{(1)})$ and  $\hat\pi_{(i)} = h(S_{(i)}) - h(S_{(i-1)}), \ i=2,\ldots, (N^2-N)/2$.  If $\hat \pi^\top o\le \hat z$, the current solution $(\hat o,\hat z)$ is optimal; otherwise, generate a valid inequality $\hat \pi^\top o \le z$.

{\color{black}{\section{Risk-Adjustable DRCC with Left-Hand Side Uncertainty}
\label{sec:lhs}
In this section, we consider  the risk-adjustable DRCC model \eqref{eq:adjust-drcc} with LHS uncertainty ($T(\xi) = \xi^\top \bar A + \bar b^\top, \ q(\xi) = q$). The DRCC \eqref{eq:drcc} is rewritten as 
\begin{equation}
\inf_{f\in\mathcal D} \mathbb P_f(A(x)\xi \ge b(x)) \ge 1-\alpha.
\end{equation}
Assume  that $A(x)\neq 0$ for all $x\in\mathcal X$.
Introducing auxiliary variables $r\in\mathbb R,\ s\in\{0,1\}^N,\  \gamma\in\mathbb R^N$ and big-M constants $M_n^1, \ n=1,\ldots,N$, the DRCC above can be reformulated as  a conic formulation (\citet{chen2022data,xie2018distributionally})
\begin{subequations}\label{eq:lhs-bilinear}
\begin{eqnarray}
	\label{eq:lhs-bilinear-constr1}		&&		\alpha N r - \sum_{n=1}^N \gamma_n \ge \epsilon N\|A(x)\|_*\\
	&& A(x)\xi^n - b(x) + M_n^1 s_n \ge r - \gamma_n, \ n=1,\ldots,N\\
	&&  (1-s_n) M_n^1  \ge r-\gamma_n, \ n=1,\ldots,N\\
	&& s_n \in \{0,1\}, \ \gamma_n \ge 0, \ n=1,\ldots,N.
\end{eqnarray}
\end{subequations}
The problem \eqref{eq:lhs-bilinear} involves a blinear term $\alpha r$ in constraint \eqref{eq:lhs-bilinear-constr1}. Formulations with bilinear constraints can be solved by GUROBI 9.0 or higher versions. The solvers relax the bilinear constraints using linear constraints based on McCormick envelops \citep{mccormick1976computability} depending on the local bounds of variables in the bilinear terms. The bounds for $\alpha$ are given by the constraint $0\le \alpha \le \bar \alpha$.
We propose tight bounds for $r$ used in generating the McCormick envelops. According to Lemma 1 in \citet{chen2022data}, the variable $r$ can be interpreted as the $\alpha$th quantile of $\{( A(x)\xi^n-b(x))^+\}_{n=1}^N$ which are in a non-descending order. Thus, the lower and upper bounds of $r$ are given by 
$$0\le r \le \left(\max_{x\in\mathcal X, n=1,\ldots,N}\{(A(x)\xi^n - b(x))\}\right)^+.$$
Consequently, the big-M constants $M_n^1$ can be chosen as $\left(\max_{x\in\mathcal X, \ell=1,\ldots,N}\{(A(x)\xi^\ell - b(x))\}\right)^+-\min_{x\in\mathcal X}\left\{A(x)\xi^n -b(x)\right\}$ for $n=1,\ldots,N$.

%
%
%
%
%
%
%
%
When the decision $x\in\mathcal X$ is restricted to binary, $\mathcal X\subset \{0,1\}^d$, we next propose an equivalent mixed integer conic (MIC) reformulation.
\begin{theorem}\label{thm:lhs}
Assume that $A(x)\neq 0$ for all $x\in \mathcal X$.
With binary decision $x\in\mathcal X\subset \{0,1\}^d$, the risk-adjustable DRCC \eqref{eq:adjust-drcc} with LHS uncertainty is equivalent to the following MIC formulation.
\begin{subequations}\label{eq:lhs-alter-final}
	\begin{eqnarray}
		\min_{x\in\mathcal X,\pi,z,\beta,t,\alpha} && c^\top x + g(\alpha)\\
		\label{eq:lhs-alter-final-constr1}		\mbox{s.t.} && 	\epsilon N{\|A(\pi) \|_*} + \sum_{n=1}^N \beta_n \le \alpha N\\
		\label{eq:lhs-alter-final-constr2}			&& A(\pi)\xi^n - b(\pi) + z_nM^2_n + \beta_n \ge 1, \ n =1,\ldots,N\\
		\label{eq:lhs-alter-final-constr3}			&& (1-z_n)M^2_n + \beta_n \ge 1, \ n =1,\ldots,N\\
		&& t \ge 0, \ \beta_n \ge 0,\ z_n\in \{0,1\}, \ n = 1,\ldots,N\\
		\label{eq:lhs-alter-final-mccormick}		&& \pi_k \ge t + (x_k - 1)M_0, \ \pi_k \le x_k M_0, \ \pi_k \le t, \ \pi_k \ge 0, \ k=1,\ldots, d\\
		\label{eq:lhs-alter-final-constr_last}	&& \alpha\in (0, \bar{\alpha}],
	\end{eqnarray}
\end{subequations}
where $M_n^2, \ n=1,\ldots,N$ and $M_0$ are big-M constants.
Furthermore, 
\begin{equation}\label{eq:lhs-valid}
	\sum_{n=1}^N z_n \le \alpha N\end{equation}
is a valid inequality.
\end{theorem}

\begin{proof}{Proof of Theorem \ref{thm:lhs}:}
According to Proposition \ref{prop:dist}, the DRCC is equivalent to 
\begin{equation}\label{eq:drcc-alt}
	\max\left\{j\in [0,N] \bigg| \sum_{n=1}^j \text{dist}(\xi^{n}, \bar S(x)) \le \epsilon \right\} \le \alpha N.\end{equation}
The largest number of elements on the left-hand side in constraint \eqref{eq:drcc-alt} corresponds to the optimal value of the (always feasible) linear program
\begin{subequations}\label{eq:drcc-al-primal}
	\begin{eqnarray}
		\max_s && \sum_{n=1}^N s_n\\
		\mbox{s.t.}&&  \sum_{n=1}^N s_n \text{dist}(\xi^n,\bar S(x)) \le \epsilon N\\
		&& 0\le s_n \le 1, \ n=1,\ldots,N.
	\end{eqnarray}
\end{subequations}
Following strong duality of linear programs, the optimal value of \eqref{eq:drcc-al-primal} coincides with that of its dual problem
\begin{subequations}\label{eq:dual-alter}
	\begin{eqnarray}
		\min_{\beta,t} && \epsilon Nt + \sum_{n=1}^N \beta_n\\
		\mbox{s.t.} && t\text{dist}(\xi^n, \bar S(x)) + \beta_n \ge 1, \ n=1,\ldots,N\\
		&& t\ge 0, \ \beta_n \ge 0, \ n=1,\ldots,N.
	\end{eqnarray}
\end{subequations}
By Lemma A.1 in \citet{chen2022data}, $\text{dist}(\xi^n,\bar S(x)) = (A(x)\xi - b(x))^+/\|A(x)\|_*$, where the convention that $0/0 = 0$ is adopted. Applying variable transformation $t \leftarrow t/\|A(x)\|_*$ and following \eqref{eq:drcc-alt} to impose the optimal value to be no more than $\alpha N$, DRCC is reformulated as the following constraints
\begin{subequations}\label{eq:dual-alter-var-trans}
	\begin{eqnarray}
		\label{eq:dual-alter-var-trans_constr1}	 && \epsilon Nt{\|A(x) \|_*} + \sum_{n=1}^N \beta_n \le \alpha N\\
		\label{eq:dual-alter-var-trans_constr2}		&& t{(A(x)\xi^n - b(x))^+} + \beta_n \ge 1, \ n=1,\ldots,N\\
		&& t\ge 0, \ \beta_n \ge 0, \ n=1,\ldots,N.
	\end{eqnarray}
\end{subequations}
The constants in \eqref{eq:dual-alter-var-trans}  have a complicating feature, the maximum operator in the second  constraints \eqref{eq:dual-alter-var-trans_constr2}, which evaluates takes the positive part of $A(x)\xi^n - b(x)$. To eliminate the maximum operator, for each $n=1,\ldots,N$, we introduce a binay variable $z_n\in \{ 0,1\}$ and express the corresponding constraint in \eqref{eq:dual-alter-var-trans_constr2} via two disjunctive  inequalities
\begin{equation}\label{eq:dual-alter-var-trans_bigM}
	t(A(x)\xi^n - b(x)) + z_nM^2_n + \beta_n \ge 1, \text{ and } (1-z_n)M^2_n + \beta_n \ge 1, \ n=1,\ldots,N.
\end{equation}
Unfortunately, the resulting constraints above and the first constraint in \eqref{eq:dual-alter-var-trans} still involve bilinear terms in $x$ and the dual variable $t$.
By exploiting the binary constraint on $x$, constraints \eqref{eq:dual-alter-var-trans} can be linearized using McCormick inequalities \eqref{eq:lhs-alter-final-mccormick} and the MIC reformulation \eqref{eq:lhs-alter-final} follows.

Next, we argue that constraint  \eqref{eq:lhs-valid} is  a valid inequality to \eqref{eq:lhs-alter-final}. To this end, take a solution $(\hat x,\hat \pi, \hat z,\hat \beta, \hat t)$ satisfying   \eqref{eq:lhs-alter-final-constr1}	-\eqref{eq:lhs-alter-final-constr_last}. We define $\bar z\in \{0,1\}^N$ such that $\bar z_n=1$ if and only if $A(x)\xi^n - b(x) <0$ for all $n=1,\ldots,N$ and claim that $(\hat x,\hat \pi, \bar z,\hat \beta, \hat t)$ satisfies \eqref{eq:lhs-alter-final-constr1}	-\eqref{eq:lhs-alter-final-constr_last} and \eqref{eq:lhs-valid}. First, we observe that $(\hat x,\bar z, \hat \beta)$ satisfies constraints \eqref{eq:lhs-alter-final-constr2}-\eqref{eq:lhs-alter-final-constr3}, which equivalently can be rewritten as
$$\min\{A(\pi)\xi^n - b(\pi) + z_nM^2_n, (1-z_n)M^2_n \} \ge 1-\beta_n, \ n=1,\ldots,N.$$
Given $\bar z$ and sufficiently large $M^2_n$  such that $A(\hat \pi)\xi^n - b(\hat \pi) + M^2_n$, we have	$\min\{A(\hat \pi)\xi^n - b(\hat \pi) + \bar z_nM_n^2, (1-\bar z_n)M_n^2 \}  = \left(A(\hat \pi)\xi^n - b(\hat \pi)\right)^+ \ge \min\{A(\hat \pi)\xi^n - b(\hat \pi) + \hat z_nM_n^2, (1-\hat z_n)M_n^2 \} $. This implies that   $(\hat x,\bar z, \hat \beta)$ satisfies constraints \eqref{eq:lhs-alter-final-constr2}-\eqref{eq:lhs-alter-final-constr3}. Now, it remains to show that $\bar z$ satisfies \eqref{eq:lhs-valid}.   To see this, we will show that $\bar z_n \le \hat \beta_n$ which then implies that $\sum_{n=1}^N \bar z_n \le \sum_{n=1}^N \hat \beta_n \le \alpha N.$ When $\bar z_n = 1$, $\min\{A(\hat \pi)\xi^n - b(\hat \pi) + \bar z_nM_n^2, (1-\bar z_n)M_n^2 \} = 0 \ge 1-\beta_n$, which implies that $\beta_n \ge 1 =\bar z_n$. When $\bar z_n = 0$, $\beta_n \ge 0 = \bar z_n$. Thus, we conclude that  \eqref{eq:lhs-valid} is  a valid inequality to \eqref{eq:lhs-alter-final}.
\end{proof}
\begin{remark}
To obtain big-M constants, consider an optimal solution $(\pi^*,t^*,\beta^*)$  to problem 	\eqref{eq:lhs-alter-final}.
We should have $t^*>0$.  To see this,  assuming that $t^*=0$, then all $\beta^*_n = 1$ for $n=1,\ldots,N$, $\pi^* = 0$ and consequently constraint \eqref{eq:lhs-alter-final-constr2} is violated. The positivity of $t^*$ implies that at least one of constraints (26b) holds tight. Thus, $t^* = \max_{n: \ A(x)\xi -b(x)>0} \{(1-\beta_n)/\left[A(x)\xi -b(x)\right]\} \le  \max_{n: \ A(x)\xi -b(x)>0} \{1/\left[A(x)\xi -b(x)\right]\}$ since $\beta_n\ge 0, \ n=1,\ldots,N $. An upper bound $\bar t$ of $t$  can be obtained as follows: for every $n=1,\ldots,N$, solve $t_n = \min_{x\in\mathcal X} \{(A(x)\xi^n -b(x))^+\}$ and let $\bar t = \max_{n: \ t_n >0}1/t_n$ if there exists at least $n$ such that $t_n >0$. The big-M constant $M_0$ can be set as $\bar t$.
The big-M constants associated with \eqref{eq:lhs-alter-final-constr2} and \eqref{eq:lhs-alter-final-constr3}  can be set as $M_n^2 = 1- \min_{x\in\mathcal X, t\ge 0}\left\{(A(x)\xi^n - b(x))t: \ \epsilon Nt\|A(x)\| \le \bar \alpha N \right\}$
%
\end{remark}}}


\section{Computational Study}
\label{sec:comp}

In the computational study, we demonstrate the computational effectiveness of the proposed mixed integer programming formulations (with both discrete and continuous distributions) on instances of a DRCC counterpart of the transportation problem with random demand \citep{luedtke2010integer,elcci2018chance}. 
For continuous distributions, 
we also compute instances of demand response management using building load where the decisions are pure binary to compare the alternative MILP (which can be found in Appendix \ref{apx:milp}) proposed in \citet{zhang2022building} and our proposed mixed 0-1 conic reformulation.
In Section \ref{sec:comp_setup}, we describe the instance setup (i) for the RHS uncertainty on the transportation problem and the demand response management problem and (ii) for the LHS uncertainty on the portfolio optimization problem. There are mainly two parts of results: (1) the computational performance (with CPU time, optimality gap, etc) of the risk-adjustable DRCC models with RHS and LHS uncertainty
in Sections \ref{sec:cpu_gaps} and \ref{sec:cpu_gaps_LHS}, respectively, and (2) the solution details given by the models in Section \ref{sec:comp_sol}. In particular, Section \ref{sec:cpu_gaps} demonstrates the computational efficacy of the proposed mixed integer formulations and valid inequalities. Section \ref{sec:comp_sol} shows that the risk-adjustable DRCC following the finite distribution assumption \ref{assump:finite} provides the highest objective values compared to  the risk-adjustable DRCC under the continuous distribution assumption \ref{assump:continuous} and the stochastic chance-constrained counterpart (which is presented in Appendix  \ref{apx:cc}).




\subsection{Computational Setup}
\label{sec:comp_setup}

\noindent\paragraph{Transportation problem:}
There are $I$ suppliers  and $D$ customers. The suppliers have limited capacity $M_i, \ i=1,\ldots, I$. There occurs a transportation cost $c_{ij}$ for shipping one unit from supplier $i$ to customer $j$. The customer demands $\tilde \xi_j , \ j=1,\ldots,D$ are random. Let $f_j$ denote the distribution of $\tilde \xi_j$ and $\mathcal D_j$ be the Wasserstein ambiguity set regarding the distribution $f_j$. 
With a penalty cost $p$ of risk tolerance $\alpha_j$ for every customer $j$, the risk-adjustable DRCC transportation problem is formulated as follows.
\begin{equation}
\min_{x\in\mathcal X,\alpha} \left\{\sum_{i=1}^I\sum_{j=1}^D c_{ij}x_{ij} +   p\sum_{j=1}^D \alpha_j: \inf_{f_j\in\mathcal D_j}\mathbb P_{f_j}(\sum_{i=1}^I x_{ij} \ge \tilde{\xi}_j) \ge 1-\alpha_j, \ 0\le \alpha_j\le \bar{\alpha}, \ j=1,\ldots,D\right\},
\end{equation}
where $\mathcal X:= \{x\in\mathbb R_+^{I\times D}: \ \sum_{j=1}^D x_{ij}\le M_i, \ i=1,\ldots,I\}$. 
Following \citet{elcci2018chance}, the risk threshold's upper bound $\bar\alpha$ is set to 0.3 and $p$ is set to $10^6$. To break symmetry, a random perturbation is added to the penalty cost $p$ following a uniform distribution on the interval $[0,100]$ for every $\alpha_j, \ j=1,\ldots,D$. The radius of the Wasserstein ball $\mathcal D_j$ is 0.05.
For other parameters (i.e., $c_{ij}, M_i$ and samples of $\tilde\xi_j$), we use the data sets with $I=40$ suppliers in \citet{luedtke2010integer} with equal probabilities for all samples.

\noindent\paragraph{Building load control problem:} There is an aggregate HVAC (i.e., heating, ventilation, and air conditioning) load of $n$ buildings to absorb random local solar photovoltaic (PV) generation $\tilde P^\text{PV}_t$ over $T$ time periods throughout the day. 
Let $\mathcal D_t$ be the Wasserstein ambiguity regarding the distribution $f$ of PV generation $\tilde P^\text{PV}_t$ during period $t$.
For each time period $t$, we solve the following risk-adjustable DRCC formulation for deciding   the room temperature $x_{t,\ell}$  and  HVAC ON/OFF decision $u_{t,\ell}$ of building $\ell$. 
\begin{equation}\label{eq:blc}
\min_{(x,u)\in\mathcal X_t,\alpha_t}\left\{ c_\text{sys}\sum_{\ell=1}^n |x_{t,\ell}-x_\text{ref}| + c_\text{switch}\sum_{\ell=1}^n u_{t,\ell} + p\alpha_t: \ \inf_{f\in\mathcal D_t}\mathbb P_f \left(\sum_{\ell=1}^n P_\ell u_{t,\ell}\ge \tilde P_t^\text{PV} \right)\ge 1-\alpha_t, \ 0\le\alpha_t\le\bar\alpha \right\},
\end{equation}
where $\mathcal X_t = \left\{(x_t,u_{t})\in\mathbb R^n\times\{0,1\}^n: \  x_{t,\ell} = A_\ell x_{t-1,\ell} + B_\ell u_{t,\ell} + G_\ell v_\ell, \ x_\text{min}\le x_{t,\ell}\le x_\text{max},  \ \ell = 1,\ldots,n \right\}$.
The objective minimizes the cost of (1) the user's discomfort (indicated by the room temperature deviation from the set-point $x_\text{ref}$), (2) switching cycles, and (3) risk violation of PV tracking. The DRCC  ensures that PV generation is absorbed by the HCAC fleet with probability $1-\alpha_t$. The indoor temperature $x_t$ and binary ON/OFF decision $u_t$ need to satisfy the constraints of temperature comfort band and thermal dynamics in the feasible set $\mathcal X_t$. The parameters $A_\ell,B_\ell,G_\ell$ are obtained from the building's thermal resistances, thermal capacity, and cooling capacity. Parameter $v_\ell$ is a given system disturbance. 
We use all the parameters and data following \citet{zhang2022building}. 
In particular, the radius of the Wasserstein ball $\mathcal D_t$ is 0.02. 
To solve the ON/OFF decisions for a planning horizon of $T=53$ periods throughout the day, one needs to sequentially solve $53$ problems in the form of
\eqref{eq:blc}, one for each period.

\noindent\paragraph{Portfolio optimization:}

Following \citet{xie2020bicriteria}, we consider a portfolio optimization problem where there are $K$ assets with random return $\tilde \xi_1,\ldots,\tilde \xi_K$. The target return is $\omega$. The investor aims to decide the investment $x_k$ into asset $k, \ k=1,\ldots,K$ with minimum cost while exceeding the target return with a high probability $1-\alpha$. The risk-adjustable DRCC portfolio optimization problem is given as following.
\begin{equation}
\min_{x\ge 0,\alpha}\left\{ \sum_{k=1}^Kc_k x_k + p\alpha: \ \inf_{f\in\mathcal D}\mathbb P_f \left(\sum_{k=1}^K\tilde \xi_k x_k \ge \omega \right) \ge 1-\alpha, \ 0\le \alpha \le \bar \alpha  \right\}
\end{equation}

We set $K\in\{30,50\}$ and $\omega=15$ and choose the cost coefficients $c_1,\ldots,c_{K}$ uniformly at random from $\{1,\ldots,100\}$. Each asset return $\tilde{\xi}_k$ is governed by a uniform distribution on $[0.8,1.5]$.
We use the $2$-norm Wasserstein ambiguity set and set Wasserstein radius to $\epsilon = 0.05$. 

The computations are conducted on 
a Windows 10 Pro machine with Intel(R) Core(TM) i7-8700 CPU3.20 GHz and 16 GB memory.
All models and algorithms are implemented
in Python 3.7.6 using Gurobi 10.0.1. 
The Gurobi default settings are used for optimizing all integer formulations except for the mixed integer conic formulation \eqref{eq:mip-continuous}, for which the Gurobi parameter \texttt{MIPFocus} is set to 3.
When implementing the branch-and-cut algorithm, we add the violated extended polymatroid inequalities using Gurobi \texttt{callback} class by \texttt{Model.cbLazy()} for integer solutions. For all the nodes in the branch-and-bound tree, we generate violated cuts at each node as long as any exists. The optimality gap tolerance is default as $10^{-4}$. The time limit  is set to 1800 seconds for computing the transportation problem instances, 100 seconds for solving the building load control problem in one period, and 3600 seconds for the portfolio optimization problem.

\subsection{CPU and Optimality Gaps with RHS Uncertainty}
\label{sec:cpu_gaps}

Under the finite distribution assumption \ref{assump:finite}, we solve the MILP \eqref{eq:mip_discrete} with and without valid inequalities in Proposition \ref{prop:finite_strengthening}. 
Under the continuous distribution assumption \ref{assump:continuous}, the mixed 0-1 conic formulation \eqref{eq:mip-continuous} can be rewritten as a mixed 0-1 second-order cone programming (SOCP) formulation if constraint  \eqref{eq:mip-continuous_conic_general} is replaced by constraints \eqref{eq:mip-continuous_conic_general_soc}. We solve the mixed 0-1 SOCP reformulation with and without valid the extended polymatroid inequalities. With only binary decisions, we also compare the mixed 0-1 SOCP reformulation with the alternative MILP reformulation in Appendix \ref{apx:milp}. Our valid inequalities significantly reduce the solution time of directly solving the mixed integer models in Gurobi.
%
%
The details are presented as follows.

\subsubsection{Finite distributions}

We first optimize transportation problem instances with the finite distribution model using the 
MILP reformulation with and without strengthening techniques proposed in Proposition \ref{prop:finite_strengthening}. 
%
Table \ref{tab:discrete_cpu} presents the CPU time (in seconds), 
``\textbf{Opt. Gap}'' as the optimality gap, and  `\textbf{`Node}'' as the total number of branching nodes.
The CPU time includes the preprocessing time $t_\text{BS}$ for calculating the violation risk $\alpha_n$ corresponding to every sample $\xi^n, \ n=1,\ldots,N$ using the bisection search method, and the time $t_\text{MILP}$ for solving the MILP reformulation \eqref{eq:mip_discrete} using Gurobi. 
In Table \ref{tab:discrete_cpu}, we solve the transportation problem with $J\in \{100,200\}$ customer demands with $N = \{50,100,200,1000,2000,3000\}$ samples. For each $(J,N)$ setting, five instances are solved. Table \ref{tab:discrete_cpu} presents the average CPU times, the average optimality gaps, and the average number of branching nodes.  Details of each instance can be found in Appendix \ref{apx:finite_cput}.

In Table \ref{tab:discrete_cpu}, with valid inequalities proposed in Proposition \ref{prop:finite_strengthening}, all the instances are solved optimally at the root node within the time limit (thus optimality gap is zero and omitted). Whereas, if being solved without the valid inequalities, the instances of more samples $(N\ge 1000)$ cannot be solved within the 18000-second time limit and ends with an  optimality gap up to 5.47\%. For larger-sized problems, solving the MILP  \eqref{eq:mip_discrete} with valid inequalities is much faster than solving the MILP  \eqref{eq:mip_discrete} directly due to the strength of the strengthened star inequality \eqref{eq:mip_discrete_constr_star}. With the valid inequalities, most of the CPU time spends on preprocessing ($t_\text{BS}$). 

\begin{table}[htbp]
\centering
\caption{Comparison of CPU time (in seconds) and optimality gaps of finite distributions}
\resizebox{.7\textwidth}{!}{
\begin{tabular}{@{\extracolsep{4pt}}ccrrrrrrrrr@{}}
\hline
\multicolumn{1}{c}{\multirow{2}[4]{*}{Demand}} & \multicolumn{1}{c}{\multirow{2}[4]{*}{$N$}} & \multicolumn{4}{c}{MILP + Valid Ineq.} &       & \multicolumn{4}{c}{MILP} \\
\cline{3-6}\cline{7-11}          &       & \multicolumn{1}{c}{$t_\text{BS}$} & \multicolumn{1}{c}{$t_\text{MILP}$} & \multicolumn{1}{c}{Time} & \multicolumn{1}{c}{Node} & \multicolumn{1}{c}{$t_\text{BS}$} & \multicolumn{1}{c}{$t_\text{MILP}$} & \multicolumn{1}{c}{Time} & \multicolumn{1}{c}{Opt. Gap} & \multicolumn{1}{c}{Node} \\
\hline
100   & 50    & 0.00  & 0.04  & 0.04  & 1     & 0.01  & 0.31  & 0.31  & N/A   & 1 \\
100   & 100   & 0.02  & 0.07  & 0.09  & 1     & 0.02  & 1.02  & 1.04  & N/A   & 158 \\
100   & 200   & 0.07  & 0.09  & 0.15  & 1     & 0.06  & 6.45  & 6.52  & N/A   & 2018 \\
100   & 1000  & 1.26  & 0.33  & 1.60  & 1     & {1.21} & {LIMIT} & {LIMIT} & {0.09\%} & {33454} \\
100   & 2000  & 4.67  & 0.72  & 5.39  & 1     & {4.64} & {LIMIT} & {LIMIT} & {0.65\%} & {9022} \\
200   & 2000  & 9.18  & 1.60  & 10.79 & 1     & {9.39} & {LIMIT} & {LIMIT} & {2.58\%} & {6678} \\
200   & 3000  & 20.37 & 2.38  & 22.75 & 1     & {20.59} & {LIMIT} & {LIMIT} & {5.47\%} & {6083} \\
\hline
\end{tabular}
}
\label{tab:discrete_cpu}%
\end{table}%

\subsubsection{Continuous distributions}
\label{sec:comp_continuous}
We first focus on the computational performance of solving the building load control problem with the binary decision. We use the proposed mixed 0-1 SOCP formulation (``\textbf{MISOCP}'') and the alternative MILP formulation (``\textbf{MILP-Binary}'') from \citet{zhang2022building}. 
In particular, the MISOCP is obtained by replacing the rotated conic constraint \eqref{eq:mip-continuous_conic} with  the SOC constraints \eqref{eq:mip-continuous_conic_general_soc1}-\eqref{eq:mip-continuous_conic_general_soc2}. 
The left-hand side function $h(o)$ of \eqref{eq:mip-continuous_conic_general_soc1} is submodular and thus the extended polymatroid inequalities \eqref{eq:extended_polymatroid_ineq} is added in a branch-and-cut (``\textbf{B\&C}'') algorithm when being violated.  

Table \ref{tab:continuous_milp_socp} reports, for each instance, the total CPU time of solving the building load control problem \eqref{eq:blc} for all 53 periods. If for any period, the problem cannot be solved within the time limit, we report ``\textbf{\#LIMIT}'' as the number of periods which cannot be solved, and ``\textbf{Avg. Gap}'' as their average optimality gap. 
Owing to its stronger relaxations and fewer variables, the MISOCP (B\&C) quickly solves all the instances, with an average of only 1.2 seconds per instance.
The optimality gaps are all zeros and thus not reported in the table. In contrast, MILP-Binary fails to be solved within the 100-second time limit for each period, with an average of 17 periods not solved to optimal. 




\begin{table}[htbp]
\centering
\caption{Comparison of CPU time (in seconds) and optimality gaps of continuous distributions with binary variables}

\resizebox{.7\textwidth}{!}{
\begin{tabular}{@{\extracolsep{4pt}}cccccccc@{}}
%
%
	\hline
	& \multicolumn{4}{c}{MILP-Binary}      & \multicolumn{3}{c}{MISOCP (B\&C)} \\
	\cline{2-5} \cline{6-8}    Instance & \multicolumn{1}{l}{Time} & \multicolumn{1}{l}{\#LIMIT} & \multicolumn{1}{l}{Avg. Node} & \multicolumn{1}{l}{Avg. Gap} & \multicolumn{1}{l}{Time} & \multicolumn{1}{l}{\#LIMIT} & \multicolumn{1}{l}{Avg. Node} \\
	\hline
	1     & 2359.28 & 21    & 314880 & 0.34\% & 1.26  & 0     & 1 \\
	2     & 2124.01 & 18    & 282644 & 0.40\% & 1.32  & 0     & 11 \\
	3     & 2165.93 & 19    & 205147 & 0.59\% & 1.30  & 0     & 1 \\
	4     & 2293.98 & 19    & 286217 & 0.62\% & 1.50  & 0     & 14 \\
	5     & 1827.85 & 14    & 202231 & 0.55\% & 1.22  & 0     & 1 \\
	6     & 2206.03 & 18    & 246280 & 0.63\% & 1.20  & 0     & 1 \\
	7     & 2190.73 & 20    & 262080 & 0.62\% & 1.23  & 0     & 18 \\
	8     & 1893.50 & 15    & 223771 & 0.33\% & 1.00  & 0     & 1 \\
	9     & 1720.97 & 13    & 211767 & 0.63\% & 1.04  & 0     & 1 \\
	10    & 1899.43 & 16    & 220385 & 0.47\% & 1.41  & 0     & 60 \\
	\hline
	avg.  & 2068.17 & 17    & 245540 & 0.52\% & 1.25  & 0     & 11 \\
	\hline
\end{tabular}%
}
\label{tab:continuous_milp_socp}%
\end{table}%

Next, we compare the branch-and-cut algorithm using the extended polymatroid inequalities (in column ``\textbf{B\&C}'') with directly solving the MISOCP \eqref{eq:mip-continuous} (in column ``\textbf{No Cuts}'') on the transportation problem instances. If any instance cannot be solved within the 1800-second time limit, we report the average optimality gap and the number of unsolved instances in parentheses. 
In Table \ref{tab:cputime_continuous_socp}, the branch-and-cut algorithm solves the MISOCP faster than directly solving it in Gurobi.
\begin{table}[htbp]
\centering
\caption{Comparison of CPU time (in seconds) and optimality gaps of continuous distributions using MISOCP}
\resizebox{.7\textwidth}{!}{
\begin{tabular}{@{\extracolsep{4pt}}clcrrrrrr@{}}
		\hline
		\multirow{2}[4]{*}{Demand} & \multirow{2}[4]{*}{$N$} & \multirow{2}[4]{*}{Instance} & \multicolumn{3}{c}{No Cuts} & \multicolumn{3}{c}{B\&C} \\
		\cline{4-6}     \cline{7-9}    
		&         &       & \multicolumn{1}{c}{Time} & \multicolumn{1}{c}{Opt. Gap} & \multicolumn{1}{c}{Node} & \multicolumn{1}{c}{Time} & \multicolumn{1}{c}{Opt. Gap} & \multicolumn{1}{c}{Node} \\
		
		\hline
		100   & 50    & a     & 92.33 & N/A   & 9186  & 7.75  & N/A   & 1 \\
		100   & 50    & b     & 66.09 & N/A   & 16896 & 10.17 & N/A   & 6 \\
		100   & 50    & c     & 80.81 & N/A   & 11245 & 10.49 & N/A   & 335 \\
		100   & 50    & d     & 38.65 & N/A   & 6526  & 15.13 & N/A   & 995 \\
		100   & 50    & e     & 67.36 & N/A   & 6565  & 8.31  & N/A   & 1 \\
		\hline
		&       & avg.  & 69.05 & N/A   & 10084 & 10.37 & N/A   & 268 \\
		\hline
		100   & 100   & a     & 64.15 & N/A   & 1     & 40.33 & N/A   & 1 \\
		100   & 100   & b     & 86.54 & N/A   & 46    & 39.79 & N/A   & 1 \\
		100   & 100   & c     & 367.28 & N/A   & 8467  & 29.57 & N/A   & 1 \\
		100   & 100   & d     & 46.37 & N/A   & 1     & 17.59 & N/A   & 1 \\
		100   & 100   & e     & 47.24 & N/A   & 1     & 29.49 & N/A   & 1 \\
		\hline
		&       & avg.  & 122.32 & N/A   & 1703  & 31.36 & N/A   & 1 \\
		\hline
		100   & 200   & a     & 1455.53 & N/A   & 1743  & 442.66 & N/A   & 1246 \\
		100   & 200   & b     & LIMIT & 0.12\% & 7211  & 434.15 & N/A   & 3 \\
		100   & 200   & c     & LIMIT & 0.32\% & 10108 & LIMIT & 0.18\% & 116520 \\
		100   & 200   & d     & LIMIT & 0.08\% & 664   & 366.55 & N/A   & 662 \\
		100   & 200   & e     & LIMIT & 0.34\% & 28139 & 619.82 & N/A   & 878 \\
		\hline
		&       & avg.  & 1731.98 & 0.21\% (4) & 9573  & 732.73 & 0.18\% (1) & 23862 \\
		\hline
	\end{tabular}%
}
\label{tab:cputime_continuous_socp}%
\end{table}%

{	\color{black}{	\subsection{CPU and Optimality Gaps with LHS Uncertainty}
	\label{sec:cpu_gaps_LHS}
	
	In this section, we solve the DRCC model \eqref{eq:adjust-drcc} with LHS uncertainty  using the conic formulations (\eqref{eq:lhs-bilinear} and \eqref{eq:lhs-alter-final}) with and without bilinear terms. A budget-type valid inequality (similar to \eqref{eq:lhs-valid}) over the binary variables is used to strengthen formulation \eqref{eq:lhs-alter-final} according to \citet{ho2023strong}.
	The portfolio optimization problem is solved with $K\in\{20,30,50\}$ assets and $N \in \{30,50\}$ samples. For each $(K,N)$ pair, we solve five instances.
	Table \ref{tab:cpu_lhs} presents the average of CPU time (in seconds), optimality gap, and total number of branching nodes over the five instances for each setting. In the column ``\textbf{Opt. Gap},'' the percentage gaps reported are averaged over the number (shown in the parentheses) of instances (out of five) which cannot be solved within the time limit of one hour. In Table \ref{tab:cpu_lhs}, MIC \eqref{eq:lhs-alter-final} is solved faster for most $(K,N)$ settings except for $(K,N)=(30,30)$. We also note that for the instance which cannot be solved within the time limit using MIC \eqref{eq:lhs-alter-final}, the optimality gaps are relatively larger than those of formulation \eqref{eq:lhs-bilinear}.
	\begin{table}[htbp]
		\centering
		\caption{Comparison of CPU time (in seconds) and optimality gaps between MIC and bilinear programs}
		\resizebox{.7\textwidth}{!}{
			\begin{tabular}{@{\extracolsep{4pt}}ccrrrrrr@{}}
				\hline
				\multirow{2}[2]{*}{K} & \multirow{2}[2]{*}{N} & \multicolumn{3}{c}{MIC \eqref{eq:lhs-alter-final}} & \multicolumn{3}{c}{Bilinear \eqref{eq:lhs-bilinear}} \\
				\cline{3-5}     \cline{6-8}    
				&       & \multicolumn{1}{l}{Time} & \multicolumn{1}{l}{Opt. Gap} & \multicolumn{1}{l}{Node} & \multicolumn{1}{l}{Time} & \multicolumn{1}{l}{Opt. Gap} & \multicolumn{1}{l}{Node} \\
				\hline
				20    & 30    & 10.22 & N/A   & 53655 & 17.76 & N/A   & 94614 \\
				20    & 50    & 184.36 & N/A   & 816809 & 543.78 & N/A   & 1067457 \\
				30    & 30    & 286.66 & N/A   & 1571725 & 151.45 & N/A   & 624048 \\
				30    & 50    & 2434.58 & 11.31\% (2) & 8636320 & 3128.93 & 5.68\% (4) & 5494919 \\
				50    & 30    & 1684.81 & 4.56\% (2) & 10433123 & 1822.76 & 3.68\% (2) & 3810056 \\
				50    & 50    & 795.12 & N/A   & 8023699 & 2929.25 & 7.27\% (3) & 4821505 \\
				\hline
			\end{tabular}%
		}
		\label{tab:cpu_lhs}%
	\end{table}%
	}}
	
	\subsection{Solution Details of Models with Finite and Continuous Distributions }
	\label{sec:comp_sol}
	In this section, we focus on the solution details of the transportation problem, which are obtained by solving the risk-adjustable DRCC models (assuming finite distributions (``\textbf{Finite}'') and continuous distributions (``\textbf{Continuous}'')), as well as the risk-adjustable stochastic chance-constrained model (``\textbf{Stochastic}''). The detailed formulation of the stochastic chance-constrained model is available in Appendix \ref{apx:cc}.
	In Section \ref{sec:cpu_gaps}, we observe that 
	the branch-and-cut algorithm does not scale as efficiently as the MILP \eqref{eq:mip_discrete}  assuming finite distributions, particularly when the sample size increases. In this section, the solution details  suggest that the MISOCP model assuming continuous distributions can be effectively approximated by the MILP model \eqref{eq:mip_discrete} for larger sample sizes.
	The details are presented below.

	\subsubsection{Optimal objective Values}
	
	We compare the optimal objective values obtained from solving the three models: Finite, Continuous, and Stochastic. In Table \ref{tab:obj_continuous_vs_discrete},  the relative difference (in columns ``\textbf{Diff.}'') is calculated as the relative gap with the Finite model. The positive relative differences of the Continuous models are as indicated by Proposition \ref{prop:obj}.
	All the relative differences  for both Continuous and Stochastic models are positive, which indicates the conservatism of the Finite model compared to the other two models. Furthermore, the differences between the Continuous and the Finite models decrease as the sample size grows larger. For instance, with a sample size $N=50$, the average  difference  between the Finite and Continuous models is 7.5\%, which reduces to 1.7\% when $N=200$. This observation implies that solving the Finite model as a conservative approximation of the Continuous model becomes more suitable when the sample size is large and the  MISOCP for the Continuous model is time-consuming to solve. 
	\begin{table}[htbp]
\centering
\caption{Comparison of objective costs}
\resizebox{.65\textwidth}{!}{
	\begin{tabular}{@{\extracolsep{4pt}}cccrrrrr@{}}
			\hline
			\multicolumn{1}{c}{\multirow{2}[4]{*}{Demand}} & \multicolumn{1}{c}{\multirow{2}[4]{*}{$N$}} & \multicolumn{1}{c}{\multirow{2}[4]{*}{Instance}} & \multicolumn{1}{c}{Finite} & \multicolumn{2}{c}{Continuous} & \multicolumn{2}{c}{Stochastic} \\
			\cline{4-4} \cline{5-6} \cline{7-8}          
			&       &       & \multicolumn{1}{c}{Obj.} & \multicolumn{1}{c}{Obj.} & \multicolumn{1}{c}{Diff.} & \multicolumn{1}{c}{Obj.} & \multicolumn{1}{c}{Diff.} \\
			\hline
			100   & 50    & a     & 37891812 & 34882415 & 7.9\% & 35877509 & 5.3\% \\
			100   & 50    & b     & 39600119 & 36583112 & 7.6\% & 37580027 & 5.1\% \\
			100   & 50    & c     & 40591537 & 37591717 & 7.4\% & 38591438 & 4.9\% \\
			100   & 50    & d     & 39992224 & 36977377 & 7.5\% & 37972322 & 5.1\% \\
			100   & 50    & e     & 41872481 & 38851630 & 7.2\% & 39851708 & 4.8\% \\
			\hline
			&       & avg.  & 39989635 & 36977250 & 7.5\% & 37974601 & 5.0\% \\
			\hline
			100   & 100   & a     & 36300456 & 34797106 & 4.1\% & 35236038 & 2.9\% \\
			100   & 100   & b     & 38370855 & 36931801 & 3.8\% & 37356067 & 2.6\% \\
			100   & 100   & c     & 39353292 & 37849666 & 3.8\% & 38297139 & 2.7\% \\
			100   & 100   & d     & 38786542 & 37271588 & 3.9\% & 37715295 & 2.8\% \\
			100   & 100   & e     & 40741978 & 39244431 & 3.7\% & 39710167 & 2.5\% \\
			\hline
			&       & avg.  & 38710625 & 37218919 & 3.9\% & 37662941 & 2.7\% \\
			\hline
			100   & 200   & a     & 35767904 & 35150000 & 1.7\% & 35202312 & 1.6\% \\
			100   & 200   & b     & 37895369 & 37270088 & 1.7\% & 37318954 & 1.5\% \\
			100   & 200   & c     & 38919823 & 38266132 & 1.7\% & 38322395 & 1.5\% \\
			100   & 200   & d     & 38302041 & 37643158 & 1.7\% & 37668983 & 1.7\% \\
			100   & 200   & e     & 40113429 & 39463734 & 1.6\% & 39516533 & 1.5\% \\
			\hline
			&       & avg.  & 38199713 & 37558622 & 1.7\% & 37605836 & 1.6\% \\
			\hline
		\end{tabular}%
	}
	\label{tab:obj_continuous_vs_discrete}
\end{table}%

\subsubsection{Allowed risk tolerance}
\begin{figure}[!htb]
	\minipage{0.32\textwidth}
	\includegraphics[width=\linewidth]{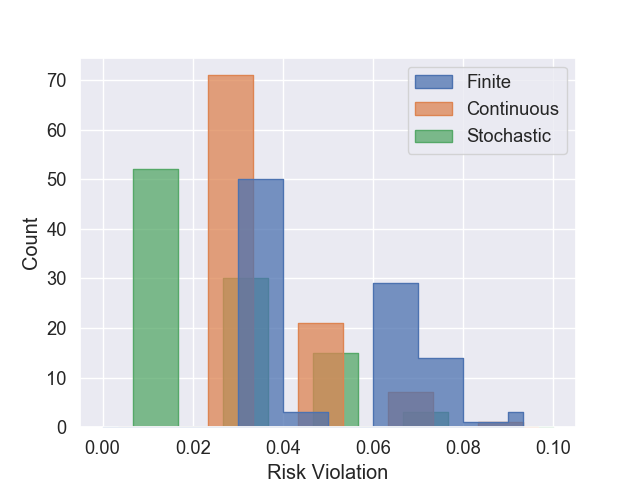}
	\caption{ $N=50$}\label{fig:awesome_image1}
	\endminipage\hfill
	\minipage{0.32\textwidth}
	\includegraphics[width=\linewidth]{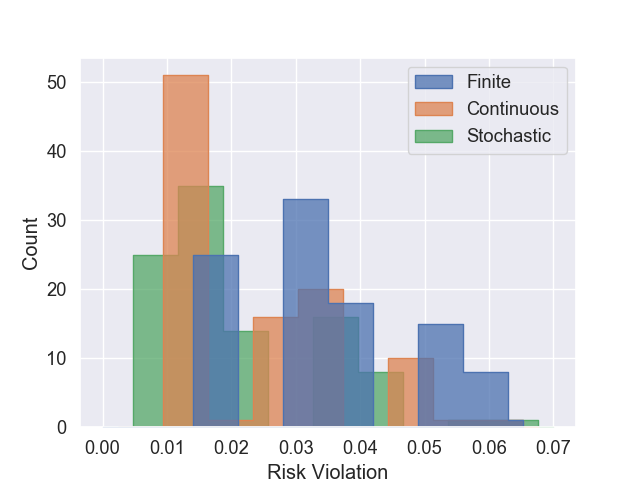}
	\caption{$N=100$}\label{fig:awesome_image2}
	\endminipage\hfill
	\minipage{0.32\textwidth}%
	\includegraphics[width=\linewidth]{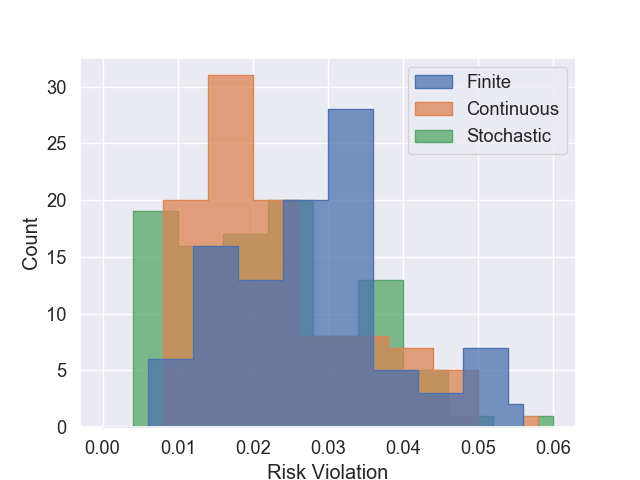}
	\caption{$N=200$}\label{fig:awesome_image3}
	\endminipage
\end{figure}

In this section, we look into the risk tolerance allowed by solving the three models. Recall that the transportation problem imposes a chance constraint for each demand location and with $D=100$ customers, there are 100 allowed risk tolerances $\alpha_j, \ j=1,\ldots,100$. Figures \ref{fig:awesome_image1}-\ref{fig:awesome_image3} show the distributions of the risk tolerances obtained by solving the Finite, Continuous, and Stochastic models with sample size $N=\{50,100,200\}$. The stochastic model assigns $\alpha$'s to smaller values than the two DRCC models. Additionally, as the sample size increases, there is more  overlap between the distributions obtained from solving the Finite and Continuous models, which supports approximating the Continuous model with the Finite model when the sample size is large. 

{	\color{black}{	\subsection{Impact of Wasserstein Ball Radius}
		
		In this section, we explore the impact of Wasserstein ball radius $\epsilon$ on the computational time and solutions using the transportation problem instances. 
		Table \ref{tab:impact_radii_bound} presents the computational times and optimal objective values of solving the Finite model with various Wasserstein ball radii $\epsilon$ and upper bounds $\bar\alpha$ of allowed risk tolerance. For a fixed radius, increasing the risk tolerance bound $\bar\alpha$ initially from 0.03 to 0.05 raises the objective value and the objective remains unchanged from 0.05 to 0.30. The objective also increases with larger radius $\epsilon$ as more ambiguity is considered. 			
		\begin{table}[htbp]
			\centering
			\caption{Comparison of CPU time (in seconds) and objectives}
			\resizebox{.45\textwidth}{!}{
				\begin{tabular}{rrrrrr}
					\hline
					\multicolumn{1}{c}{$\epsilon$} & \multicolumn{1}{c}{ $\bar\alpha$} & \multicolumn{1}{c}{Obj} & \multicolumn{1}{c}{$t_\text{BS}$} & \multicolumn{1}{c}{$t_\text{MILP}$} & \multicolumn{1}{c}{Time} \\
					\hline
					0.01  & 0.03  & 75338651.10 & 0.26  & 0.35  & 0.61 \\
					0.01  & 0.05  & 75305108.16 & 0.66  & 0.51  & 1.17 \\
					0.01  & 0.30  & 75305108.16 & 20.08 & 2.66  & 22.74 \\
					\hline
					0.05  & 0.03  & 75720278.48 & 0.26  & 0.30  & 0.56 \\
					0.05  & 0.05  & 75681030.16 & 0.68  & 0.42  & 1.10 \\
					0.05  & 0.30  & 75681030.16 & 20.77 & 2.51  & 23.28 \\
					\hline
					0.10  & 0.03  & 76005476.44 & 0.25  & 0.28  & 0.53 \\
					0.10  & 0.05  & 75960989.41 & 0.72  & 0.41  & 1.13 \\
					0.10  & 0.30  & 75960989.41 & 20.65 & 2.41  & 23.06 \\
					\hline
				\end{tabular}%
			}
			\label{tab:impact_radii_bound}%
		\end{table}%
		Figure \ref{fig:awesome_epsilon} displays the distributions and medians of the allowed risk tolerances obtained for radii $\epsilon\in\{0.01,0.05,0.10\}$ with an upper risk tolerance bound $\bar\alpha = 0.30$. 
		The risk-adjustable model assigns $\alpha$ to larger values with larger radii to better trade-off between the system cost and risk violation.
		\begin{figure}[!htb]
			\centering
			\includegraphics[width=.4\linewidth]{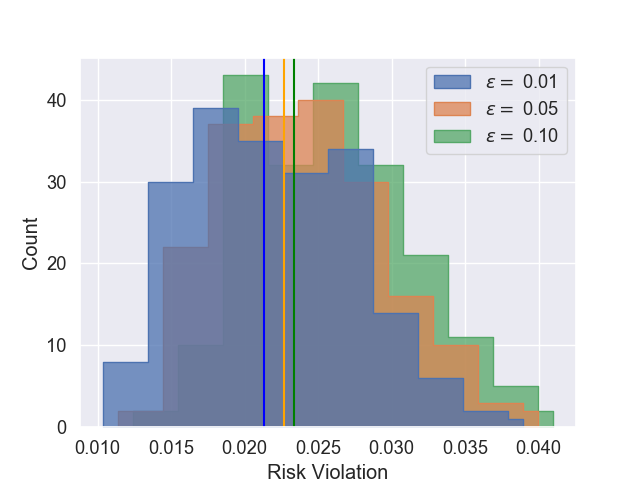}
			\caption{ $\bar\alpha=0.30$}
			\label{fig:awesome_epsilon}
		\end{figure}
}}

\section{Conclusions}
\label{sec:conclusions}

In this paper, we investigated distributionally robust individual chance-constrained problems with a data-driven Wasserstein ambiguity set, where the uncertainty  affects either the RHS or the LHS. The risk tolerance is treated as a decision variable. The goal of the risk-adjustable DRCC is to trade-off between system costs and risk violation costs via penalizing the risk tolerance in the objective function.
For the RHS uncertainty,  we provided a MILP reformulation of the risk-adjustable DRCC problem with finite distributions and a MISOCP reformulation for the continuous distribution case. 
For the LHS uncertainty, we mixed integer conic reformulations with binary decisions.
Moreover, valid inequalities are derived for all the reformulations. Via extensive numerical studies, we demonstrated that our valid inequalities accelerate solving the risk-adjustable DRCC models when compared to optimization solvers. 
Although the MISOCP reformulation does not scale well with larger  sample size, the MILP reformulation can be used as an approximation of the MISOCP reformulation.




%
	\bibliographystyle{informs2014}
	\bibliography{YilingRef}
	
	\newpage
	\begin{APPENDICES}
		
		\section{Dominance of Risk Tolerance}
		\label{sec:dominance}
		The chance constrained programming literature \citep[see, e.g.,][]{prekopa1990dual,dentcheva2000concavity,ruszczynski2002probabilistic,prekopa2003probabilistic} defines the concept of \emph{non-dominated points}, or the so-called $p$-efficient points, where $p$ refers to $1-\alpha$ in this paper. 
		
		
		\begin{definition}[$p$-efficient point] \citep{prekopa2003probabilistic,dentcheva2000concavity}
			Let $p\in(0,1)$. A point $v\in\mathbb R^m$ is a $p$-efficient point of the probability distribution $f$, $\mathbb P_f(v) \ge p$ and there is no $w\le v, \ w\neq v$ such that $\mathbb P_f(w) \ge p$.
		\end{definition}
		The concept can be extended to the DRO variant in the following definition.
		\begin{definition}[Distributionally Robust $p$-efficient point]
			Let $p\in(0,1)$. A point $v\in\mathbb R^m$ is a distributionally robust $p$-efficient point of the ambiguity set $\mathcal D$, $\inf_{f\in\mathcal D}\mathbb P_f(v) \ge p$ and there is no $w\le v, \ w\neq v$ such that $\inf_{f\in\mathcal D}\mathbb P_f(w) \ge p$.
		\end{definition}
		
		In the case of individual chance constraint with an uncertain RHS, under the empirical distribution of $\{\xi^n\}_{n=1}^N$, the $(1-\alpha)$-efficient point is the $(1-\alpha)$-quantile (or $(1-\alpha)$-VaR) of the empirical distribution. The distributionally robust $(1-\alpha)$-efficient point coincides with the worst-case VaR $t_\alpha$ (obtained by assuming either the finite distribution or the continuous distribution), which is greater than the $(1-\alpha)$-quantile of the reference distribution in the Wasserstein ball $\mathcal D$. Similar to the $(1-\alpha)$-quantile, the worst-case VaR is nonincreasing in the risk tolerance, which is formally stated in the following proposition.
		\begin{proposition}\label{prop:VaR}
			Given $0\le\alpha_1< \alpha_2\le\bar\alpha$, let $t_{\alpha_1}$ and $t_{\alpha_2}$ be the worst-case VaRs associated with $\alpha_1$  and $\alpha_2$, respectively. Then, $t_{\alpha_1} \ge t_{\alpha_2}$.
		\end{proposition}
		The proof can be easily derived based on the water-filling interpretation in Section \ref{sec:prelim} and  is omitted for brevity.
		
		
		
		\section{Stochastic Chance Constrained Problem:  MILP Reformulation}
		\label{apx:cc}
		Let $N^\prime = \lceil\bar\alpha N\rceil$. \begin{subequations}\label{eq:mip_sp}
			\begin{eqnarray}
				\min_{x\in\mathcal X, y} && c^\top x + \sum_{n=1}^{N^\prime-1}\Delta_n y_n + \Delta_{N^\prime}  \\
				\label{eq:sp_constr1}\mbox{s.t.} && Tx \ge \xi^{N^\prime} + \sum_{n=1}^{N^\prime - 1}(\xi^n - \xi^{n+1}) y_n\\
				\label{eq:sp_y_order}&& y_{n+1} \ge y_{n}, \ n=1,\ldots,N^\prime -2\\
				\label{eq:sp_alpha}&& \frac{1}{N} \left( N^\prime-1- \sum_{n=1}^{N^\prime-1}y_{n} \right) \in (0,\bar{\alpha}]\\
				&& y_n\in\{0,1\}, \ n=1,\ldots,N^\prime - 1,
			\end{eqnarray}
		\end{subequations}
		where $\Delta_n:=g(n/N) - g((n+1)/N), \ n=1,\ldots,N^\prime-1$, and $\Delta_{N^\prime} := g(N^\prime /N)$.

		\section{Alternative MILP Reformulation for Risk-adjustable DRCC with Binary Variables}
		\label{apx:milp}
		When all the decision variables are pure binary, i.e., $\mathcal X\subset \{0,1 \}^d$, \citet{zhang2022building} developed a MILP-base reformulation. The reformulation uses a binary variable $o_{jk}$ to identify the critical indices $j$ and $k$ following the similar idea as in Section \ref{sec:continuous}. 
		\begin{equation}\label{eq:soc-milp}
			\sum_{j=0}^{N-1}\sum_{k=j}^{N-1}o_{jk}\left[(\alpha N-j)(Tx-\xi^{k+1}) - \sum_{i=j+1}^k(\xi^i-\xi^{k+1}) \right]\ge N\epsilon 
		\end{equation}
		%
		%
		There are   bilinear terms $o_{jk}x_\ell, \alpha o_{jk}$ and trilinear term $o_{jk}\alpha x_\ell$ in constraint \eqref{eq:soc-milp}. When the decisions $x_\ell, \ \ell=1,\ldots,d$ are pure binary, they can all be linearized using McCormick inequalities \citep{mccormick1976computability}. The alternative MILP reformulation is as follows.
		\begin{subequations}\label{eq:milp-continuous}
			\begin{eqnarray}
				\min_{x\in\mathcal X, o,\alpha,u,w}&& c^\top x + g(\alpha)\\
				\mbox{s.t.} && \sum_{j=0}^{N-1}\sum_{k=j}^{N-1}\left[ o_{jk}\sum_{i=j+1}^{k} \left( \xi^{k+1}-\xi^i \right) + N T(\delta_{jk} - j\tau_{jk}) - \xi^{k+1}(\varepsilon_{jk} - j o_{jk}) \right] \ge N\epsilon \\
				&& \sum_{j=0}^{N-1}\sum_{k=j}^{N-1} \xi^{j}o_{jk} \ge Tx \ge \sum_{j=0}^{N-1}\sum_{k=j}^{N-1} \xi^{j+1}o_{jk}\\
				&& \sum_{j=0}^{N-1} \sum_{k=j}^{N-1} (k+1) o_{jk} \ge \alpha N \ge \sum_{j=0}^{N-1} \sum_{k=j}^{N-1} k o_{jk}\\
				\label{eq:milp-continuous_sos1}&& \sum_{j=0}^{N-1}\sum_{k=j}^{N-1} o_{jk}=1\\
				&& \alpha \in (0,\bar{\alpha}]\\
				&&  \varepsilon_{jk}\le o_{jk}, \ \varepsilon_{jk}\le \alpha, \ \varepsilon_{jk}\ge \alpha + o_{jk} -1, \ \varepsilon_{jk} \ge 0, \ 0\le j\le k\le N-1 \\
				&&  \delta_{\ell jk}\le \varepsilon_{jk}, \ \delta_{\ell jk}\le x_\ell, \ \delta_{\ell jk}\ge \varepsilon_{jk} + x_{\ell} -1, \ \delta_{\ell jk} \ge 0, \ 0\le j\le k\le N-1, \ 1 \le \ell \le d \\
				&&  \tau_{\ell jk}\le o_{jk}, \ \tau_{\ell jk}\le x_\ell, \ \tau_{\ell jk}\ge o_{jk} + x_{\ell} -1, \ \tau_{\ell jk} \ge 0, \ 0\le j\le k\le N-1, \ 1 \le \ell \le d \\
				&& o_{jk}\in \{0,1\}, \ 0\le j\le k \le N-1.
			\end{eqnarray}
		\end{subequations}
		
		\section{More results for CPU time and Optimality Gaps with Finite Distributions}
		\label{apx:finite_cput}
		\begin{table}[htbp]
			\centering
			\caption{Comparison of CPU time (in seconds) and optimality gaps of finite distributions}
			\resizebox{.6\textwidth}{!}{
				\begin{tabular}{@{\extracolsep{4pt}}cccrrrrrrrrr@{}}
					\hline
					\multicolumn{1}{c}{\multirow{2}[4]{*}{Demand}} & \multicolumn{1}{c}{\multirow{2}[4]{*}{$N$}} & \multicolumn{1}{c}{\multirow{2}[4]{*}{Instance}} & \multicolumn{4}{c}{MILP + Valid Ineq.} &       & \multicolumn{4}{c}{MILP} \\
					\cline{4-7}\cline{8-12}  &        &       & \multicolumn{1}{c}{$t_\text{BS}$} & \multicolumn{1}{c}{$t_\text{MILP}$} & \multicolumn{1}{c}{Time} & \multicolumn{1}{c}{Node} & \multicolumn{1}{c}{$t_\text{BS}$} & \multicolumn{1}{c}{$t_\text{MILP}$} & \multicolumn{1}{c}{Time} & \multicolumn{1}{c}{Opt. Gap} & \multicolumn{1}{c}{Node} \\	
					\hline
					100   & 50    & a     & 0.02  & 0.05  & 0.06  & 1     & 0.01  & 0.33  & 0.34  & N/A   & 1 \\
					100   & 50    & b     & 0.00  & 0.05  & 0.05  & 1     & 0.01  & 0.25  & 0.26  & N/A   & 1 \\
					100   & 50    & c     & 0.00  & 0.03  & 0.03  & 1     & 0.01  & 0.31  & 0.32  & N/A   & 1 \\
					100   & 50    & d     & 0.00  & 0.03  & 0.03  & 1     & 0.00  & 0.38  & 0.38  & N/A   & 1 \\
					100   & 50    & e     & 0.00  & 0.05  & 0.05  & 1     & 0.00  & 0.27  & 0.27  & N/A   & 1 \\
					\hline
					&       & avg.  & 0.00  & 0.04  & 0.04  & 1     & 0.01  & 0.31  & 0.31  & N/A   & 1 \\
					\hline
					100   & 100   & a     & 0.02  & 0.06  & 0.08  & 1     & 0.02  & 0.91  & 0.93  & N/A   & 85 \\
					100   & 100   & b     & 0.02  & 0.06  & 0.08  & 1     & 0.02  & 0.84  & 0.87  & N/A   & 1 \\
					100   & 100   & c     & 0.02  & 0.08  & 0.09  & 1     & 0.02  & 1.06  & 1.08  & N/A   & 223 \\
					100   & 100   & d     & 0.02  & 0.07  & 0.09  & 1     & 0.03  & 1.37  & 1.40  & N/A   & 480 \\
					100   & 100   & e     & 0.02  & 0.08  & 0.10  & 1     & 0.02  & 0.91  & 0.93  & N/A   & 1 \\
					\hline
					&       & avg.  & 0.02  & 0.07  & 0.09  & 1     & 0.02  & 1.02  & 1.04  & N/A   & 158 \\
					\hline
					100   & 200   & a     & 0.06  & 0.09  & 0.15  & 1     & 0.07  & 4.73  & 4.81  & N/A   & 1420 \\
					100   & 200   & b     & 0.06  & 0.09  & 0.14  & 1     & 0.06  & 5.41  & 5.47  & N/A   & 1705 \\
					100   & 200   & c     & 0.06  & 0.08  & 0.14  & 1     & 0.06  & 6.87  & 6.94  & N/A   & 2605 \\
					100   & 200   & d     & 0.08  & 0.08  & 0.16  & 1     & 0.06  & 8.21  & 8.28  & N/A   & 2559 \\
					100   & 200   & e     & 0.06  & 0.09  & 0.16  & 1     & 0.06  & 7.04  & 7.11  & N/A   & 1800 \\
					\hline
					&       & avg.  & 0.07  & 0.09  & 0.15  & 1     & 0.06  & 6.45  & 6.52  & N/A   & 2018 \\
					\hline
					100   & 1000  & a     & 1.20  & 0.33  & 1.54  & 1     & 1.20  & LIMIT & LIMIT & 0.06\% & 55681 \\
					100   & 1000  & b     & 1.23  & 0.35  & 1.58  & 1     & 1.19  & LIMIT & LIMIT & 0.05\% & 32960 \\
					100   & 1000  & c     & 1.23  & 0.32  & 1.55  & 1     & 1.20  & LIMIT & LIMIT & 0.06\% & 25663 \\
					100   & 1000  & d     & 1.47  & 0.31  & 1.78  & 1     & 1.23  & LIMIT & LIMIT & 0.23\% & 26593 \\
					100   & 1000  & e     & 1.19  & 0.35  & 1.54  & 1     & 1.21  & LIMIT & LIMIT & 0.05\% & 26371 \\
					\hline
					&       & avg.  & 1.26  & 0.33  & 1.60  & 1     & 1.21  & LIMIT & LIMIT & 0.09\% & 33454 \\
					\hline
					100   & 2000  & a     & 4.62  & 0.74  & 5.36  & 1     & 4.69  & LIMIT & LIMIT & 0.70\% & 7529 \\
					100   & 2000  & b     & 4.70  & 0.75  & 5.45  & 1     & 4.70  & LIMIT & LIMIT & 0.59\% & 9567 \\
					100   & 2000  & c     & 4.64  & 0.69  & 5.33  & 1     & 4.58  & LIMIT & LIMIT & 0.67\% & 8824 \\
					100   & 2000  & d     & 4.72  & 0.70  & 5.42  & 1     & 4.61  & LIMIT & LIMIT & 0.64\% & 11390 \\
					100   & 2000  & e     & 4.68  & 0.70  & 5.38  & 1     & 4.62  & LIMIT & LIMIT & 0.67\% & 7799 \\
					\hline
					&       & avg.  & 4.67  & 0.72  & 5.39  & 1     & 4.64  & LIMIT & LIMIT & 0.65\% & 9022 \\
					\hline
					200   & 2000  & a     & 9.12  & 1.59  & 10.71 & 1     & 9.27  & LIMIT & LIMIT & 1.80\% & 6840 \\
					200   & 2000  & b     & 9.25  & 1.61  & 10.86 & 1     & 9.34  & LIMIT & LIMIT & 1.86\% & 6835 \\
					200   & 2000  & c     & 9.19  & 1.58  & 10.77 & 1     & 9.43  & LIMIT & LIMIT & 3.98\% & 6615 \\
					200   & 2000  & d     & 9.26  & 1.67  & 10.93 & 1     & 9.48  & LIMIT & LIMIT & 3.24\% & 6479 \\
					200   & 2000  & e     & 9.08  & 1.58  & 10.66 & 1     & 9.45  & LIMIT & LIMIT & 2.00\% & 6619 \\
					\hline
					&       & avg.  & 9.18  & 1.60  & 10.79 & 1     & 9.39  & LIMIT & LIMIT & 2.58\% & 6678 \\
					\hline
					200   & 3000  & a     & 20.20 & 2.16  & 22.36 & 1     & 20.88 & LIMIT & LIMIT & 5.09\% & 6549 \\
					200   & 3000  & b     & 20.08 & 2.58  & 22.66 & 1     & 20.22 & LIMIT & LIMIT & 5.50\% & 6510 \\
					200   & 3000  & c     & 20.45 & 2.39  & 22.84 & 1     & 20.54 & LIMIT & LIMIT & 4.81\% & 6574 \\
					200   & 3000  & d     & 20.73 & 2.36  & 23.09 & 1     & 20.81 & LIMIT & LIMIT & 6.00\% & 4232 \\
					200   & 3000  & e     & 20.39 & 2.39  & 22.78 & 1     & 20.51 & LIMIT & LIMIT & 5.96\% & 6552 \\
					\hline
					&       & avg.  & 20.37 & 2.38  & 22.75 & 1     & 20.59 & LIMIT & LIMIT & 5.47\% & 6083 \\
					
					\hline
				\end{tabular}
			}
			\label{tab:discrete_cpu_details}%
		\end{table}%

	\end{APPENDICES}
	%
	%

	
	

\end{document}